\documentclass[11pt,a4paper]{article}
\usepackage[margin=24mm,headheight=14pt]{geometry}
\usepackage[T1]{fontenc}
\usepackage[utf8]{inputenc}
\usepackage{lmodern}
\usepackage{amsmath,amssymb,amsthm,mathtools}
\usepackage{booktabs,array,longtable,tabularx,multirow}
\usepackage{microtype,needspace,placeins}
\usepackage{xcolor,listings,enumitem}
\usepackage{fancyhdr}
\usepackage{fvextra}
\usepackage[unicode,pdfencoding=auto]{hyperref}
\usepackage{bookmark,xurl}
\definecolor{linkblue}{RGB}{30,63,91}
\hypersetup{colorlinks=true,linkcolor=linkblue,citecolor=linkblue,urlcolor=linkblue,
 pdftitle={Improved Lower Bounds for Mutually Orthogonal Latin Squares (MOLS) at 33 Orders},
 pdfsubject={Explicit constructions of mutually orthogonal Latin squares}}
\allowdisplaybreaks[2]
\setlist{nosep,leftmargin=1.6em}
\newtheorem{theorem}{Theorem}[section]

\theoremstyle{definition}
\newcommand{\F}{\mathbb F}
\newcommand{\Z}{\mathbb Z}
\newcommand{\OA}{\operatorname{OA}}
\newcommand{\im}{\operatorname{im}}
\newcommand{\enc}{\operatorname{enc}}

\newcommand{\xor}{\mathbin{\oplus}}
\newcommand{\code}[1]{\texttt{\detokenize{#1}}}

\newcolumntype{Y}{>{\raggedright\arraybackslash}X}
\title{Improved Lower Bounds for\\
Mutually Orthogonal Latin Squares (MOLS) at 33 Orders}
\author{
    Christopher D. Rosin \\
    https://constructive.codes \\
    christopher.rosin@gmail.com
}
\date{}
\begin{document}
\maketitle
\begin{abstract}
A \textit{Latin square} of order $n$ is an $n\times n$ array over symbols $\{0,1,\ldots, n-1\}$, such that the same symbol does not appear twice in any row or any column.  Two Latin squares of order $n$ are \textit{orthogonal} if, when superimposed, they yield each of the $n^2$ ordered pairs of symbols exactly once.  A set of Latin squares are \textit{mutually orthogonal Latin squares (MOLS)} if each pair of distinct squares is orthogonal.  Let $N(n)$ be the maximum size of a set of MOLS of order $n$.  For most $n$ which are not prime powers, there is a wide gap between known upper and lower bounds on $N(n)$.  This has led to interest in constructing new larger sets of MOLS for explicit values of $n$, to incrementally improve the lower bounds on $N(n)$.  This paper provides MOLS that improve the lower bound on $N(n)$ for 33 values of $n$ in the range from 22 through 494.  The constructions were found by GPT-6 Astra, and use established techniques including difference matrices, quasi-difference matrices, orthogonal arrays, and projective planes.  A highlight is 4 MOLS of order $n=22$; combined with pre-existing results, this establishes $N(n) \geq 4$ for all $n>10$.  
\end{abstract}

\section{Introduction}

Mutually orthogonal Latin squares (MOLS) are a well-studied type of combinatorial design \cite{Handbook}.  A \textit{Latin square} of order $n$ is an $n\times n$ array over symbols $\{0,1,\ldots, n-1\}$, such that the same symbol does not appear twice in any row or any column.  Two Latin squares of order $n$ are \textit{orthogonal} if, when superimposed, they yield each of the $n^2$ ordered pairs of symbols exactly once (for a small example, see the end of ``A small example: two MOLS of order 4'' in Section~\ref{sec:affine-arrays}).  A set of Latin squares are \textit{mutually orthogonal Latin squares (MOLS)} if each pair of distinct squares is orthogonal.  Let $N(n)$ be the maximum size of a set of MOLS of order $n$. 

The central problem in the study of MOLS is to find $N(n)$ for all $n$, or at least limit the range of possible values with good lower and upper bounds.  $N(n) \leq n-1$ for all $n\geq 2$, and for prime powers $q=p^k$ for prime $p$ and $k \geq 1$, constructions using finite fields and projective planes achieve $N(q)=q-1$ \cite{Bose}.  Outside of prime powers though, there remains a wide gap between upper and lower bounds.  For $n=10$ and also for non-prime-power $n \geq 14$ with $n\equiv 1$ or $2$ $($mod $ 4)$ and with $n$ not a sum of two squares, we have $N(n) \leq n-4$, with only slightly better upper bounds for larger $n$ satisfying these conditions \cite{BruckRyser,Shrikhande,Bruck,Lam}.  For other non-prime-power $n \geq 12$, $N(n) \leq n-1$ remains the only known upper bound.  

For lower bounds, $N(n) \geq 2$ for all $n>6$ \cite{boseshrikhandeparker1960}, and $N(n) \geq 3$ for all $n>10$ \cite{Wallis84,Todorov}.  In this paper, we improve this to $N(n) \geq 4$ for all $n>10$ by constructing for the first time 4 MOLS of order $n=22$, completing the last open case.  

Somewhat larger lower bounds are known for specific $n$.  A 2024 compilation \cite{mols500} of lower bounds for each $n<500$ noted improvements in lower bounds for 9 orders $n<120$ since the publication of the 2007 \textit{Handbook of Combinatorial Designs} \cite{Handbook}.  It appears that no new improvements in lower bounds for $n<500$ have appeared in the literature since the 2024 compilation.  MOLS construction has been included as a task in recent AI benchmarks \cite{molsai1,molsai2}, and the results in this paper show how AI approaches can indeed construct new sets of MOLS.  Compared to the 2024 compilation, this paper provides improved lower bounds for 14 orders $n<120$, and for another 19 orders in $120 \leq n < 500$.  

Order 10 illustrates some of the difficulty involved in pinning down $N(n)$.  We know that $2 \leq N(10) \leq 6$, but substantial ongoing computer search efforts \cite{mols10a,mols10b,mols10c} have not yet narrowed this.  Even for 3 MOLS of order $n=10$, and even with well-developed search techniques, the search space is still too large to search completely.

Improvements to lower bounds have generally proceeded via highly structured objects that can be succinctly described. The constructions below use techniques along these lines: difference matrices (Section~\ref{sec:dm}) and quasi-difference matrices (Sec.~\ref{sec:qdm}) with additional symmetry constraints, lifting coarse orthogonal arrays (Sec.~\ref{sec:affine-arrays}), and selecting points from a projective plane (Sec.~\ref{sec:geometric}).  

\section{Results}

The new lower bounds are listed in Table~\ref{tab:results}, with a comparison to the previous lower bounds from the 2024 compilation \cite{mols500}.  Each set of MOLS is available as a text file in the repository \texttt{https://github.com/Constructive-Codes/MOLS}.  It is straightforward to check that each set of MOLS has the appropriate size for MOLS of order $n$, uses only symbols in $\{0,1,\ldots, n-1\}$, has no repeated symbol in any row or any column of any square, and meets the requirement that any pair of distinct squares when superimposed contains each ordered pair of symbols exactly once.  The repository provides Python code that performs the verification.

\begin{theorem}\label{thm:main}
For each order $n$ and associated lower bound in the Constructed column in Table~\ref{tab:results}, $N(n)$ is at least as large as the given lower bound.
\end{theorem}

The 4 MOLS of order $n=22$, combined with pre-existing lower bounds \cite{mols500} and \cite[Theorem 5.1]{wilson}, give:
\begin{theorem}\label{thm:bound4}
$N(n)\geq 4$ for all $n>10$.
\end{theorem}

In addition, the 6 MOLS at order $n=74$ combined with pre-existing lower bounds \cite{mols500} and \cite[Theorem 5.1]{wilson}, give:
\begin{theorem}\label{thm:bound6}
$N(n)\geq 6$ for all $n>68$.
\end{theorem}

The 4 MOLS of order $n=22$ enable construction of a strongly regular graph with parameters (484,126,40,30); the existence of this was previously open \cite{Brouwertable}.  Via the standard construction \cite{brouwersrgbook}, the 484 vertices correspond to the $22 \times 22 = 484$ cells in a Latin square; two vertices are adjacent when their cells share a row, a column, or a symbol in one of the 4 squares.  The adjacency matrix for this SRG is provided in the repository as well.

The constructions of the MOLS are described in the sections below.

\section{AI Usage}

The large language model OpenAI GPT-6 Astra \cite{gpt6astra} found all the MOLS constructions.  Both ChatGPT Pro, and the OpenAI Codex harness \cite{openaicodex}, were used with GPT-6 Astra.  Prompts are shown in Appendix~\ref{app:prompts}.  

The descriptions of the constructions in sections below were written by GPT-6 Astra, with human author review and light editing.  The accompanying code implementing these constructions was written by GPT-6 Astra, and is provided in this paper's repository  \texttt{https://github.com/Constructive-Codes/MOLS}.  Each set of MOLS was verified by code written by the author, with this code also provided in the repository.

For additional verification of the construction descriptions below, a separate large language model Claude Opus 5.5 \cite{opus55} was given the construction sections and data appendices below \textit{without} the accompanying code, and succeeded in constructing MOLS giving the same results as in Table~\ref{tab:results}; the Opus 5.5 constructed MOLS were separately verified.

\begin{center}\small
\begin{longtable}{rrrrl}
\caption{New lower bounds on $N(n)$, the size of the largest set of MOLS of order $n$, from constructions in this paper.  Each new lower bound is compared to the previously-compiled lower bound from 2024 \cite{mols500}.}\label{tab:results}\\
\toprule $n$ & 2024 & Constructed & Gain & Construction\\\midrule\endfirsthead
\toprule $n$ & 2024 & Constructed & Gain & Construction\\\midrule\endhead
22 & 3 & 4 & 1 & Quasi-difference matrix\\
45 & 6 & 7 & 1 & Difference matrix\\
50 & 6 & 7 & 1 & Lift coarse orthogonal array\\
54 & 8 & 9 & 1 & Lift coarse orthogonal array\\
60 & 5 & 6 & 1 & Quasi-difference matrix\\
72 & 7 & 8 & 1 & Difference matrix\\
74 & 5 & 6 & 1 & Quasi-difference matrix\\
78 & 6 & 7 & 1 & Relative difference matrix\\
80 & 9 & 10 & 1 & Difference matrix\\
85 & 6 & 7 & 1 & Difference matrix\\
90 & 6 & 8 & 2 & Quasi-difference matrix\\
92 & 6 & 7 & 1 & Relative difference matrix\\
96 & 10 & 11 & 1 & Difference matrix\\
98 & 6 & 7 & 1 & Lift coarse orthogonal array\\
144 & 10 & 14 & 4 & Difference matrix\\
162 & 6 & 9 & 3 & Lift coarse orthogonal array\\
175 & 6 & 7 & 1 & Projective plane\\
183 & 6 & 7 & 1 & Projective plane\\
192 & 7 & 9 & 2 & Difference matrix\\
207 & 8 & 9 & 1 & Projective plane\\
225 & 8 & 9 & 1 & Difference matrix\\
288 & 15 & 17 & 2 & Difference matrix\\
324 & 8 & 11 & 3 & Difference matrix\\
328 & 7 & 9 & 2 & Projective plane\\
335 & 6 & 9 & 3 & Projective plane\\
427 & 7 & 12 & 5 & Projective plane\\
430 & 6 & 11 & 5 & Projective plane\\
442 & 7 & 11 & 4 & Projective plane\\
456 & 7 & 11 & 4 & Projective plane\\
469 & 8 & 11 & 3 & Projective plane\\
482 & 6 & 7 & 1 & Projective plane\\
486 & 8 & 9 & 1 & Lift coarse orthogonal array\\
494 & 6 & 11 & 5 & Projective plane\\
\bottomrule
\end{longtable}
\end{center}

\section{Common representations and development algorithms}\label{sec:framework}
\subsection{Symbols, groups, and field arithmetic}
Unless an explicit range states otherwise, indices start at zero.
An additive product group
$G=\Z_{r_0}\times\cdots\times\Z_{r_{d-1}}$ uses the integer label
\[
 \enc(x_0,\ldots,x_{d-1})
 =x_0+r_0x_1+r_0r_1x_2+\cdots.
\]
Decode the label, perform addition or subtraction coordinatewise modulo
the corresponding radices, and encode again. For the cyclic group $\Z_m$,
addition and subtraction are integer arithmetic modulo $m$. For $\F_p^d$ the
labels use base-$p$ encoding with coordinate zero least significant; in particular,
binary vector addition is XOR. Matrices over $\F_p$ act on column vectors.
An explicitly displayed tuple is read in its displayed coordinate order.

For an extension field, the same digits encode polynomial coefficients.
Multiply polynomials and reduce by the defining polynomial specified for
the construction. Keep in mind that integer labels in an extension field are not residues
modulo its cardinality. Unless a construction gives a different order,
enumerate integer labels increasingly, with the leftmost named loop outermost.
Any consistent relabeling preserves the mathematical result, but these
conventions also make the numerical output reproducible.

\subsection{Orthogonal arrays as lists of grid cells}\label{sec:oa}
An \emph{index-one orthogonal array of strength two}, written
$\OA(n^2,w,n,2)$, is a list of $n^2$ records with $w$ coordinates in an
$n$-symbol alphabet. In any two coordinates, each ordered symbol pair
must occur exactly once. This is simply a pair-uniqueness condition on a
table, so the orthogonal array is equivalent to a set of MOLS. For every record $(a_0,a_1,\ldots,a_{w-1})$, set the Latin square
\begin{equation}\label{eq:extract}
 L_j(a_0,a_1)=a_{j+2},\qquad 0\le j<w-2.
\end{equation}
The pair $(a_0,a_1)$ addresses each cell once. The coordinate pairs $(0,j+2)$ and
$(1,j+2)$ establish the row and column properties; pairs $(i+2,j+2)$
establish orthogonality. Thus $w$ array coordinates give $w-2$ squares.
The implementation is a loop that fills initially empty matrices.

\subsection{Difference matrices as short translation recipes}\label{sec:dm}
Let $G$ be an additive abelian group of size $m$. An \emph{index-one
difference matrix} (DM) with $w$ columns is an $m\times w$ matrix $D$
such that, for each $i<j$,
\begin{equation}\label{eq:dm}
 \{D_{r,j}-D_{r,i}:0\le r<m\}=G
 \quad\text{with each element occurring once.}
\end{equation}
Subtract $D_{r,0}$ from every entry of row $r$, giving
$\bar D_{r,j}=D_{r,j}-D_{r,0}$. Define
\begin{equation}\label{eq:dm-squares}
 L_j(r,g)=\bar D_{r,j}+g,\qquad 1\le j<w,\quad g\in G.
\end{equation}
Each row of $L_j$ is a translate of $G$. Each column of $L_j$ is a
permutation of $G$ because, for $j>0$, column $j$ of $\bar D$ lists
every group element once.
For two specified outputs $x,y$, the difference $y-x$ determines $r$
uniquely by~\eqref{eq:dm}, and then $g=x-\bar D_{r,i}$ is unique.
Hence the algorithm gives $w-1$ MOLS. Equivalently,
$(r,\bar D_{r,0}+g,\ldots,\bar D_{r,w-1}+g)$ is an OA with $w+1$
coordinates. This explains the one-coordinate difference between the
DM and OA counts~\cite{abelcheng1994,greevenbroek2018}.
Abel~\cite{abel2015} used this development to obtain five MOLS of
order 60 from a $(60,6,1)$ difference matrix.

\subsection{Quasi-difference matrices and filling a hole}\label{sec:qdm}
Difference-matrix development uses a group as its symbol set. To construct
squares of an order larger than the group, we can adjoin new symbols and
leave a smaller part of the orthogonal array to be filled separately.
A quasi-difference matrix describes this extension.

Let $G$ be a finite abelian group, written additively, and let $m=|G|$.
Adjoin $h$ distinct symbols
$I=\{\infty_0,\ldots,\infty_{h-1}\}$, none of which belongs to $G$.
These are labels, not group elements: adding $g\in G$ to a row changes
its entries in $G$ and leaves every $\infty_a$ fixed.
The final alphabet will be $G\cup I$, of size $n=m+h$.

We use the following form of a \emph{quasi-difference matrix} (QDM):
an $(m+2h)\times w$ matrix $Q$ over $G\cup I$ such that
\begin{enumerate}
\item each row contains at most one entry from $I$;
\item each symbol $\infty_a$ occurs exactly once in each column;
\item for any two columns $i<j$, the differences $Q_{r,j}-Q_{r,i}$,
      taken over rows $r$ with both entries in $G$, contain every
      element of $G$ exactly once.
\end{enumerate}
For a fixed column pair, the first two conditions exclude exactly $2h$
rows from the difference list: $h$ have a new symbol in the first column
and $h$ have one in the second. Thus the list in the third condition has
$m$ entries, as required.

\paragraph{Development and the missing pairs.}
For each row $r$ of $Q$ and each $g\in G$, form a row whose entry in
column $i$ is
\[
 q^{(r,g)}_i=\begin{cases}
 Q_{r,i}+g,&Q_{r,i}\in G,\\
 \infty_a,&Q_{r,i}=\infty_a.
 \end{cases}
\]
These $m(m+2h)$ rows cover every ordered pair of symbols except pairs
in $I\times I$. In the first two coordinates, these missing pairs
correspond to an empty $h\times h$ region in each extracted square,
called a \emph{hole}. To fill it, append the rows of an
$\OA(h^2,w,h,2)$, replacing its
$h$ symbols by $\infty_0,\ldots,\infty_{h-1}$. The resulting array has
\[
 m(m+2h)+h^2=(m+h)^2=n^2
\]
rows. We check that each ordered pair occurs exactly once in every
two columns $i,j$.

If the prescribed pair is $(x,y)\in G\times G$, translation preserves
its difference. The third QDM condition therefore gives exactly one
seed row $r$ with $Q_{r,j}-Q_{r,i}=y-x$. In that row,
$g=x-Q_{r,i}$ is the unique translation giving $x$ in column $i$;
it also gives $y$ in column $j$.

If the prescribed pair is $(\infty_a,y)$ with $y\in G$, the second
condition gives exactly one seed row $r$ with $Q_{r,i}=\infty_a$.
Its entry $Q_{r,j}$ belongs to $G$, by the first condition. The unique
translation $g=y-Q_{r,j}$ gives $y$ in column $j$ while leaving
$\infty_a$ in column $i$. The same argument applies when the finite
symbol is in column $i$ instead.

Finally, no translated seed row contains two symbols from $I$. Every
prescribed pair $(\infty_a,\infty_b)$ therefore occurs only in the
hole filling, where the OA property gives it exactly once. This proves
that the completed array is an $\OA(n^2,w,n,2)$ and hence gives $w-2$
MOLS by~\eqref{eq:extract}.

\paragraph{A small example.}
Take $G=\Z_3$, $h=1$, and $w=3$. The matrix
\[
 Q=\begin{pmatrix}
 \infty&0&2\\
 0&\infty&1\\
 0&2&\infty\\
 0&0&0\\
 0&1&2
 \end{pmatrix}
\]
satisfies the QDM conditions: each column contains $\infty$ once,
and each finite column-pair difference list is a permutation of
$0,1,2$. For example, in columns $0,1$ the three finite differences
are $2,0,1$. Developing the five rows by $g=0,1,2$ gives fifteen
rows; appending $(\infty,\infty,\infty)$ completes the array.
The requested pair $(\infty,1)$ in columns $0,1$ selects the first
seed and translation $g=1$, giving $(\infty,1,0)$.
With $\infty$ labelled 3, extraction gives the Latin square
\[
 L=\begin{pmatrix}
 0&2&3&1\\
 3&1&0&2\\
 1&3&2&0\\
 2&0&1&3
 \end{pmatrix}.
\]
In particular, the appended row fills cell $(3,3)$ with symbol 3.

\paragraph{Filling a hole with a finite field.}
For $h=1$, the single all-infinity row works for any number of columns.
For a prime-power $h>1$, choose distinct nonzero elements
$a_1,\ldots,a_{w-2}$ of $\F_h$; this requires $w\le h+1$.
For each $(u,v)\in\F_h^2$, form
\begin{equation}\label{eq:hole}
 (u,v,u+a_1v,\ldots,u+a_{w-2}v).
\end{equation}
Any two specified coordinates determine $u$ and $v$ uniquely.
For the first two coordinates this is immediate. Given $u$ and
$u+a_iv=x$, solve $v=(x-u)/a_i$; given $v$ and $u+a_iv=x$, solve
$u=x-a_iv$. For two coordinates of the form
$u+a_iv=x$ and $u+a_jv=y$, subtract the equations to obtain
\[
 v=\frac{y-x}{a_j-a_i},\qquad u=x-a_iv.
\]
The divisions are valid because $a_i\ne0$ and $a_j\ne a_i$.
Thus every ordered pair of field symbols occurs once in each
coordinate pair. Replacing field label $b$ by $\infty_b$ gives
the required hole filling.

The constructions below use all six nonzero labels when $h=7$ and all
eight when $h=9$. At $h=9$ the field is $\F_3[z]/(z^2+1)$,
with label $a+3b$ for $a+bz$. In the final squares, $\infty_b$
has integer label $m+b$.

Quasi-difference matrices were used by Todorov~\cite{todorov1989}
to construct four MOLS of order 20 and by Abel and
Todorov~\cite{abeltodorov1993} for orders 20, 30, 38, and 44.
Further constructions and hole fillings appear in Abel~\cite{abel2008}.

\section{Affine partitions: a common lifting algorithm}
Many of the constructions below describe each symbol by two components:
a symbol from a small \emph{coarse} alphabet and a vector from a \emph{fine}
space over a prime field. A short array specifies the coarse components.
Each of its rows is then expanded into several rows by varying an input
vector. Linear maps and fixed offsets determine the fine components of
these new rows. The resulting array is much larger, but its defining data
remain small.

The central question is whether this expansion covers the required values
exactly once. For a difference matrix we must cover all differences; for
an orthogonal array we must cover all ordered pairs. In both cases, the
linear maps generate subspaces, and the offsets place translated copies
of these subspaces so that they partition the required space. If $U$ is a
linear subspace, a translated copy $c+U=\{c+u:u\in U\}$ is called an
\emph{affine subspace}, or a \emph{coset} of $U$. Two such cosets are either
equal or disjoint. This elementary fact allows us to check an entire
expanded array by checking considerably smaller lists of offsets.

\subsection{Lifting a coarse difference array}\label{sec:affine-lifts}
We first construct a difference matrix over $S\times V$, where $S$ is an
abelian group of size $s$ and $V=\F_p^d$. Its order will be
$n=sp^d$. A matrix with $r$ columns will then give $r-1$ MOLS by the
development algorithm in Section~\ref{sec:dm}.

Choose an input space $Y=\F_p^a$, with $0\le a\le d$, and put
$\lambda=p^{d-a}$. We will replace each coarse row by $p^a$ rows, one for
each $y\in Y$, so the coarse array must have
$n/p^a=s\lambda$ rows. The ingredients have three distinct roles:
\begin{itemize}
\item An $(s\lambda)\times r$ array $H$ over $S$ supplies the coarse
      component $H_{t,i}$ of column $i$ in row $t$.
\item A linear map $M_i:Y\to V$ for each column specifies how its fine
      component changes as $y$ varies. The same map is used in every
      coarse row of that column.
\item An offset $c_{t,i}\in V$ for each entry of $H$ specifies the starting
      point of that fine component.
\end{itemize}
For each coarse row $t$ and each input $y$, form the row with entries
\begin{equation}\label{eq:affine-lift}
 D_{(t,y),i}=(H_{t,i},M_i y+c_{t,i})\in S\times V.
\end{equation}
Thus the coarse component stays fixed while the fine component runs
through an affine subspace. The expansion has
$s\lambda p^a=sp^d=n$ rows, as required.

To see what is needed for these rows to form a difference matrix, fix
two distinct columns $i,j$. Their difference in row $(t,y)$ is
\[
 \bigl(H_{t,i}-H_{t,j},\ (M_i-M_j)y+c_{t,i}-c_{t,j}\bigr).
\]
The coarse difference restricts which rows $t$ can be used. Within one
such row, the fine differences run through a translate of
$U_{ij}=\im(M_i-M_j)$. We impose the following conditions for every
column pair:
\begin{enumerate}
\item Each $h\in S$ occurs exactly $\lambda$ times among
      $H_{t,i}-H_{t,j}$. Thus $H$ is a difference matrix of index
      $\lambda$: every group difference occurs $\lambda$ times.
\item The map $M_i-M_j$ is injective. Consequently $U_{ij}$ has
      $p^a$ elements, and varying $y$ produces each element once.
\item For each fixed $h\in S$, the $\lambda$ cosets
\[
 (c_{t,i}-c_{t,j})+U_{ij},
 \qquad t\text{ with }H_{t,i}-H_{t,j}=h,
\]
      are distinct. Since each has $p^a$ elements and
      $\lambda p^a=|V|$, they partition $V$.
\end{enumerate}
All three conditions must hold for the same arrays and maps, for every
pair of columns. In particular, the offsets cannot be chosen independently
for different pairs, because each offset participates in several tests.

These conditions give a direct uniqueness proof. Suppose a difference
$(h,v)\in S\times V$ is prescribed. First restrict to the $\lambda$
coarse rows with difference $h$. Their cosets partition $V$, so exactly
one of these rows has a coset containing $v$. In that row, injectivity
of $M_i-M_j$ gives exactly one input $y$ producing $v$. Thus every group
difference occurs exactly once. Normalize the first column of $D$ to
zero and apply~\eqref{eq:dm-squares} to obtain the $r-1$ squares.

For the numerical construction, enumerate $t$ first and then $y$ in
increasing base-$p$ encoding. Encode a group element $(h,v)$ as
$\enc(h)+s\enc_p(v)$, where
$\enc_p(v_0,\ldots,v_{d-1})=\sum_bv_bp^b$.

\paragraph{Checking the partition with short labels.}
The coset condition can be checked without listing every vector in every
coset. Once $M_i-M_j$ is known to be injective, Gaussian elimination finds
a linear map
\[
 R_{ij}:V\longrightarrow\F_p^{d-a},\qquad \ker R_{ij}=U_{ij}.
\]
For example, the rows of a matrix for $R_{ij}$ can be chosen as a basis
of the linear forms that vanish on every column of $M_i-M_j$.
Two vectors $v,v'$ have the same label under $R_{ij}$ exactly when
$v-v'\in U_{ij}$, which is exactly the condition that they belong to
the same coset. Hence, for each coarse difference $h$, it is enough to
check that
\[
 R_{ij}(c_{t,i}-c_{t,j}),
 \qquad t\text{ with }H_{t,i}-H_{t,j}=h,
\]
lists every vector of $\F_p^{d-a}$ once. These labels describe the
\emph{quotient space} $V/U_{ij}$: they identify which coset contains
a vector, rather than its position within that coset.

For binary fine spaces with $a=d-1$, each label is a single bit and each
coarse difference must select one coset of each label. This is the
affine-coset mechanism of Abel and Cheng~\cite[Sections 5--7, especially
conditions 6.1--6.3]{abelcheng1994}, who obtained, among other results,
five MOLS of order 48, seven of order 40, and nine of order 80. The
formulation above uses the same partition principle for other prime
characteristics and larger quotients.

\paragraph{A small example: two MOLS of order 12.}
Take $S=\Z_3$, $V=\F_2^2$, and $Y=\F_2$, so $\lambda=2$.
Index the six coarse rows by $(q,u)\in\Z_3\times\F_2$, and use three
columns. Define
\[
 H_{(q,u),*}=(0,q,2q),\qquad
 M_0y=(0,0),\quad M_1y=(y,0),\quad M_2y=(0,y),
\]
and
\[
 c_{(q,u),0}=(0,0),\qquad
 c_{(q,u),1}=(0,u),\qquad
 c_{(q,u),2}=(u,u).
\]
The coarse differences $H_1-H_0$, $H_2-H_0$, and $H_2-H_1$ are
$q$, $2q$, and $q$, respectively. Thus any prescribed coarse difference
fixes $q$ and leaves the two choices $u=0,1$. The corresponding fine
differences, as $y$ runs through $0,1$, are as follows; all vector
arithmetic in the table is over $\F_2$.
\begin{center}
\begin{tabular}{cll}
\toprule
Column difference & $u=0$ & $u=1$\\
\midrule
$D_1-D_0$ & $\{(0,0),(1,0)\}$ & $\{(0,1),(1,1)\}$\\
$D_2-D_0$ & $\{(0,0),(0,1)\}$ & $\{(1,1),(1,0)\}$\\
$D_2-D_1$ & $\{(0,0),(1,1)\}$ & $\{(1,0),(0,1)\}$\\
\bottomrule
\end{tabular}
\end{center}
Every row of this table partitions all four vectors of $V$, and its
two-element sets also show that the difference maps are injective.
This checks all three column pairs. In integer labels, the twelve rows
are given completely by
\[
 D_0=0,\qquad D_1=q+3y+6u,\qquad
 D_2=(2q\bmod3)+3u+6((y+u)\bmod2).
\]
Enumerate $q=0,1,2$, then $u=0,1$, then $y=0,1$. For example,
$(q,u)=(0,1)$ expands to $(0,6,9)$ and $(0,9,3)$.
Developing this $12\times3$ difference matrix gives two MOLS of order
12. The larger constructions use the same steps with more columns and
larger fine spaces.

\subsection{Lifting a coarse orthogonal array}\label{sec:affine-arrays}
The second construction produces an orthogonal array directly. Its
alphabet again consists of a coarse symbol and a fine vector, but now
we must make every \emph{ordered pair of complete symbols} occur once
in every pair of columns. This requires $n^2$ output rows and will give
$w-2$ MOLS when there are $w$ columns.

Label the coarse alphabet $0,\ldots,s-1$, let the fine space be
$V=\F_p^d$, and set $n=sp^d$. Choose an input space $X=\F_p^e$ with
$0\le e\le2d$, and put $\lambda=p^{2d-e}$. A single coarse row will
produce $p^e$ output rows. For any fixed pair of coarse symbols, there
are $p^{2d}$ possible pairs of fine vectors, so we need $\lambda$
coarse rows to provide them.

Accordingly, use an $s^2\lambda$-by-$w$ array $H$ in which every ordered
pair of coarse symbols occurs $\lambda$ times in every column pair.
This is a strength-two OA of index $\lambda$. As before, attach an
offset $C_{t,i}\in V$ to each entry. The linear maps have a little more
freedom: for each column $i$ and coarse symbol $a$, choose a map
$A_{i,a}:X\to V$. Thus a column may use different maps in rows with
different coarse symbols. Within any fixed group of rows having the
same coarse symbols in columns $i,j$, the two relevant maps are fixed,
while the offsets may vary.

For every coarse row $t$ and every input $x\in X$, output the row
\begin{equation}\label{eq:affine-oa}
 \left(H_{t,i}+s\enc_p\bigl(A_{i,H_{t,i}}x+C_{t,i}\bigr)\right)_{i=0}^{w-1}.
\end{equation}
The coarse component is the remainder modulo $s$, and the remaining
part encodes the fine vector. Enumerate $t$ first and then $x$ in
increasing base-$p$ encoding. There are
$s^2\lambda p^e=s^2p^{2d}=n^2$ rows.

Fix columns $i<j$ and coarse outputs $(a,b)$, and write
\[
 T_{ij}^{ab}=\{t:(H_{t,i},H_{t,j})=(a,b)\},\qquad
 Bx=(A_{i,a}x,A_{j,b}x)\in V\times V.
\]
Here $B$ is the linear map obtained by stacking the two matrices, so it
has $2d$ rows and $e$ columns. The required checks have the same form
as in the difference-matrix construction:
\begin{enumerate}
\item $|T_{ij}^{ab}|=\lambda$, as required by the coarse OA.
\item $B$ is injective, or equivalently has column rank $e$. A single
      map $A_{i,a}$ need not be injective; the two outputs together must
      determine the input.
\item The $\lambda$ cosets
\[
 (C_{t,i},C_{t,j})+\im B,\qquad t\in T_{ij}^{ab},
\]
      are distinct and hence partition $V\times V$, since each has
      $p^e$ elements and $\lambda p^e=p^{2d}$.
\end{enumerate}
These checks are made for every column pair and every ordered pair of
coarse symbols.

For the correctness argument, prescribe two complete output symbols.
Their coarse components determine $(a,b)$ and therefore the set
$T_{ij}^{ab}$. Their fine components belong to exactly one of the
displayed cosets, determining $t$. Injectivity of $B$ then determines
$x$ uniquely. Thus every ordered pair occurs exactly once, and the
expanded array is an $\OA(n^2,w,n,2)$. Apply~\eqref{eq:extract} to obtain
the $w-2$ squares.

\paragraph{The hyperplane test.}
At orders 50, 54, and 98, the coarse alphabet is binary and
$e=2d-1$, so $\lambda=p$; their fine spaces have characteristics
$p=5,3,7$, respectively. Once $B$ has rank $2d-1$, its image has one
dimension fewer than $V\times V$. Such a subspace is called a
\emph{hyperplane}. Gaussian elimination gives a nonzero vector
$\ell\in\F_p^{2d}$ such that $\ell^TB=0$. The hyperplane is exactly
\[
 \im B=\{z\in V\times V:\ell^Tz=0\}:
\]
the image lies in this kernel, and both spaces have dimension $2d-1$.
For an offset $c=(C_{t,i},C_{t,j})$, every point in its translated image
has the same label, since $\ell^T(Bx+c)=\ell^Tc$. Consequently this
coset consists exactly of the vectors with the single label
\[
 \ell^Tz=\ell^T(C_{t,i},C_{t,j}).
\]
The $p$ possible labels describe all $p$ cosets. The partition test is
therefore reduced to checking that these offset labels, for
$t\in T_{ij}^{ab}$, are $0,1,\ldots,p-1$ in some order. The choice of
nonzero scalar multiple of $\ell$ only permutes the labels and does not
affect the test.

\paragraph{A small example: two MOLS of order 4.}
Take $s=p=2$, $d=e=1$, and $w=4$, so $V=X=\F_2$ and $\lambda=2$.
Index the eight coarse rows by $(t,u,v)\in\F_2^3$ and define
\[
 H_{(t,u,v),*}=(t,\ t+v,\ t+u+v,\ t+u),\qquad
 C_{(t,u,v),*}=(0,\ u,\ v,\ u+v).
\]
All additions in these tuples are modulo two. Let $A_{i,a}x=x$ for
every $i,a$. Each coarse row produces two output rows, with entries
$H_i+2(x+C_i)$ for $x=0,1$, where the addition inside the parentheses
is again modulo two. These formulas specify all sixteen rows.

For any column pair, $B:x\mapsto(x,x)$ has image
$\{(0,0),(1,1)\}$. Its two cosets are distinguished by
$\ell=(1,1)$, so the offset label is $C_i+C_j$.
Consider columns 0 and 1 with coarse outputs $(0,0)$. The equations
$H_0=H_1=0$ force $t=v=0$ and leave $u=0,1$. Their offsets are
$(0,0)$ and $(0,1)$; varying $x$ gives the two disjoint sets
\[
 \{(0,0),(1,1)\}\quad\hbox{and}\quad
 \{(0,1),(1,0)\},
\]
which cover all four pairs of fine outputs.

The same check for all six column pairs can be read from the following
table. Each row represents two column pairs.
\begin{center}
\begin{tabular}{cll}
\toprule
Column pairs $(i,j)$ & $H_i+H_j$ & $C_i+C_j$\\
\midrule
$(0,1),(2,3)$ & $v$ & $u$\\
$(0,2),(1,3)$ & $u+v$ & $v$\\
$(0,3),(1,2)$ & $u$ & $u+v$\\
\bottomrule
\end{tabular}
\end{center}
In each row, the last two entries are independent linear forms in
$(u,v)$. Fixing $H_i=a$ and $H_j=b$ fixes the first of these forms to
$a+b$. There are two possibilities for $(u,v)$, and the offset label
is zero in one and one in the other. Each possibility determines a
unique $t$ from $H_i=a$. Thus the coarse pair occurs twice and its two
fine cosets are distinct. This verifies every required pair and proves
an $\OA(16,4,4,2)$.

Using the first two coordinates as row and column addresses gives
\[
 L_0=\begin{pmatrix}
 0&3&1&2\\2&1&3&0\\3&0&2&1\\1&2&0&3
 \end{pmatrix},\qquad
 L_1=\begin{pmatrix}
 0&2&3&1\\3&1&0&2\\1&3&2&0\\2&0&1&3
 \end{pmatrix}.
\]

\paragraph{Common translations and subgroup development.}
For orders 98, 162, and 486, part of the input adds the same fine
translation to every column. Write the input as $(y,z)$, with
$y\in\F_p^a$, $z\in V$, and use maps of the form
\[
 A_{i,\alpha}(y,z)=M_{i,\alpha}y+z.
\]
Here $e=a+d$ and $\lambda=p^{d-a}$. At orders 162 and 486 the maps
$M_{i,\alpha}=M_i$ are independent of the coarse bit. At order 98 the
input $y$ has one coordinate over $\F_7$, and the maps also depend on
the coarse bit; its construction below specifies them explicitly.

For fixed coarse outputs $(\alpha,\beta)$, subtracting the two fine
outputs cancels $z$. The pair test therefore reduces to the following
condition: $M_{i,\alpha}-M_{j,\beta}$ is injective, and the $\lambda$ sets
\[
 \{(M_{i,\alpha}-M_{j,\beta})y+C_{t,i}-C_{t,j}:y\in\F_p^a\},
 \qquad t\in T_{ij}^{\alpha\beta},
\]
partition $V$. A prescribed fine difference determines $t,y$; either
individual fine output then determines $z$. This is the difference
partition used in~\eqref{eq:affine-lift}, now applied separately to
each ordered coarse pair.

There is also a direct difference-matrix interpretation. In these three
cases the coarse alphabet is $\Z_2$. Before adding $z$, the construction
has $4\lambda p^a=4p^d=2n$ seed rows over $\Z_2\times V$. Fix a pair
of columns and a requested group difference $(h,v)$. There are two
ordered binary pairs $(\alpha,\beta)$ with $\alpha-\beta=h$, and the
partition condition gives one seed row for each. Thus the seeds form an
\emph{index-two difference matrix}: each group difference occurs twice.
The two occurrences have different binary components in either of the
chosen columns, so their entries lie in different cosets of
$K=\{0\}\times V$. Translating each seed by the $|K|=n/2$ elements
of $K$ produces the $n^2$ OA rows.
This is the subgroup-development mechanism used by Todorov to construct
four MOLS of order 14~\cite{todorov2012}, and formulated in terms of a
group and subgroup by Abel and Li~\cite[Section 2]{abelli2015}.

\section{Small affine difference matrices}
\begin{center}
\begin{tabular}{rclrrrrr}
\toprule
$n$ & $S$ & $V$ & $a$ & $\lambda$ & rows of $H$ & $r$ & MOLS\\
\midrule
45 & $\Z_5$ & $\F_3^2$ & 1 & 3 & 15 & 8 & 7\\
72 & $\F_2^3$ & $\F_3^2$ & 1 & 3 & 24 & 9 & 8\\
80 & $\Z_5$ & $\F_2^4$ & 2 & 4 & 20 & 11 & 10\\
96 & $\Z_3$ & $\F_2^5$ & 3 & 4 & 12 & 12 & 11\\
\bottomrule
\end{tabular}
\end{center}

To label the entries, encode $(v_0,\ldots,v_{d-1})\in\F_p^d$ as
$\sum_b v_b p^b$ and encode $(h,v)\in S\times V$ as
$\operatorname{enc}(h)+s\operatorname{enc}(v)$.
At order 72 the coarse encoding is the three-bit integer, whose addition
is bitwise XOR; in the other three cases it is the usual residue modulo
$s$.  Each map is specified by the encoded columns of its $d\times a$
matrix.  All four constructions enumerate $t$ first and then $y$ in
increasing encoded order.  Appendix~\ref{app:affine-certificates} supplies
all numerical data.

\subsection{Order 45: ternary lines over a five-symbol array}
Use the $15\times8$ array $H$ in Appendix~\ref{app:affine-certificates},
an index-three difference matrix over $\Z_5$ of the type constructed by
Kuhfeld~\cite{kuhfeld2005}. The eight maps are $M_i y=y\mu_i$, where
\[
 (\operatorname{enc}(\mu_0),\ldots,\operatorname{enc}(\mu_7))
   =(1,3,6,5,4,2,8,7).
\]
Decode every entry into two ternary digits.  Use the $15\times8$ offset
table and \eqref{eq:affine-lift}.  Every coarse difference selects three
affine lines, each with three points; the displayed offsets make them
partition $\F_3^2$.  This yields a $45\times8$ difference matrix over
$\Z_5\times\F_3^2$ and seven MOLS.

\subsection{Order 72: a coarse permutation program and affine recurrence}
For a three-bit integer $x$, define
\[
 T(x)=((2x)\bmod8)\mathbin{\mathrm{XOR}}
 \begin{cases}3,&x\ge4,\\0,&x<4.\end{cases}
\]
This is a linear permutation of $\F_2^3$ with $T^7=I$.
Take the permutations, written in input order $0,\ldots,7$,
\[
\begin{split}
 f_0&=(1,6,3,0,7,5,2,4),\\
 f_1&=(2,4,3,7,1,0,6,5),\\
 f_2&=(7,6,4,3,2,1,0,5).
\end{split}
\]
For $t=0,1,2$ and $x=0,\ldots,7$, let
\[
 H_{8t+x}=\bigl(0,x,\ T^j f_t(T^{-j}x)\text{ for }j=0,\ldots,6\bigr).
\]
Coarse subtraction is XOR.  For the fine maps $M_i y=y\mu_i$, use
\[
 (\operatorname{enc}(\mu_i))_{i=0}^8=(0,1,2,3,4,6,5,7,8).
\]
The appendix gives eight seed rows $c_{x,i}\in\F_3^2$.
To generate the other sixteen rows, define
\[
 S=(0,7,1,6,4,3,5,2),\quad \tau=(2,0,1),\quad
 P(0)=0,\ P(1)=1,\ P(2+j)=2+(2j\bmod7).
\]
Whenever $c_{8t+x,i}=(u,v)$ is known, assign
\begin{equation}\label{eq:72-recurrence}
 c_{8\tau(t)+S(x),P(i)}=(u+v,\ v+\delta_i),\qquad
 \delta_i=\begin{cases}i,&i<3,\\0,&i\ge3,\end{cases}
 \quad\text{in }\F_3^2.
\end{equation}
Apply this first to layer $t=0$ to fill layer $2$, and then to layer $2$
to fill layer $1$.  Formula~\eqref{eq:affine-lift} now gives the
$72\times9$ difference matrix over $\F_2^3\times\F_3^2$.
The data satisfy all 288 affine partition conditions.

\subsection{Order 80: four-point subspaces in a sixteen-element field}
Use $\F_{16}=\F_2[X]/(X^4+X+1)$, representing a polynomial by the
integer whose bits are its coefficients.  Addition is XOR.  The four
field elements $U=\{0,1,6,7\}$ form its subfield of order four.
Identify $y=(y_0,y_1)\in\F_2^2$ with $y_0+6y_1\in U$.

The coarse array has four blocks of five rows.  Given the vectors
$A_s,B_s$ printed in the appendix, set, for $s=0,\ldots,3$,
$x,j\in\{0,\ldots,4\}$,
\[
 H_{5s+x,0}=0,\quad H_{5s+x,1+j}=A_s[(x+j)\bmod5],\quad
 H_{5s+x,6+j}=B_s[(x+j)\bmod5].
\]
The eleven field slopes are
\[
 (a_0,\ldots,a_{10})=(0,1,15,10,12,8,11,3,2,13,7).
\]
Set $M_i y=a_i(y_0+6y_1)$, or equivalently give $M_i$ the encoded
binary columns $(a_i,6a_i)$, with multiplication in $\F_{16}$.
Use the offset table in the appendix.  Each coarse difference selects
four translates of $(a_i+a_j)U$; the displayed offsets make them partition
$\F_{16}$.  The result is an $80\times11$ difference matrix over
$\Z_5\times\F_2^4$ and ten MOLS.

\subsection{Order 96: binary three-spaces over a ternary array}
Use the $12\times12$ coarse array $H$ in the appendix.  It is obtained
from a 36-row, 13-column ternary strength-two array in Sloane's
catalogue~\cite{sloane36} by taking
rows $0,3,\ldots,33$ and columns $0,\ldots,11$.
Every coarse column difference takes each value in $\Z_3$ four times.
The appendix lists twelve maps $M_i:\F_2^3\to\F_2^5$ by their three
five-bit columns, and all $12\times12$ five-bit offsets.
Apply \eqref{eq:affine-lift}.  Each difference map is injective and
each coarse difference selects four affine three-spaces, each of size
eight, partitioning $\F_2^5$.
This gives a $96\times12$ difference matrix over
$\Z_3\times\F_2^5$ and eleven MOLS.

\section{Binary coarse arrays for orders 50, 54 and 98}
Each instance uses the orthogonal-array lift~\eqref{eq:affine-oa} with $s=2$.

\subsection{Order 50}
Use $(p,d,e,w)=(5,2,3,9)$.  In this case the two fine maps for a
coordinate are equal: $A_{i,0}=A_{i,1}=A_i$.
The following three seeds determine all nine matrices:
\[
 A_0=\begin{bmatrix}4&1&0\\1&0&1\end{bmatrix},\qquad
 A_3=\begin{bmatrix}0&1&0\\2&0&1\end{bmatrix},\qquad
 A_6=\begin{bmatrix}1&0&0\\0&2&1\end{bmatrix}.
\]
Let $\rho(a,b,c)=(c,a,b)$.  For $b\in\{0,3,6\}$ and
$s\in\{0,1,2\}$, obtain $A_{b+s}$ by applying $\rho^s$ to
each row of $A_b$.  On a nine-entry vector let $R$ apply $\rho$
separately to its three consecutive blocks of three entries.
Table~\ref{tab:affine50-seeds} gives the binary rows and offsets in eight
short orbits.  An offset $(c_0,c_1)$ is packed as
$D=c_0+5c_1$, so recover it as $(D\bmod5,\lfloor D/5\rfloor)$.
For a seed of length $\ell$ beginning at row $t$, append
$(H_{t+s,*},D_{t+s,*})=(R^sH_{t,*},R^sD_{t,*})$ for
$0\leq s<\ell$, in the displayed seed order.
The two all-zero binary rows have different offsets. This gives
$20$ binary rows and $180$ offset vectors.

Equation~\eqref{eq:affine-oa} now produces
$20\cdot5^3=2500$ rows and hence seven MOLS of order 50.
For example, for coordinates $0,1$ the vector
$\ell=(0,1,4,4)^T$ annihilates the stacked matrix on the left.
For binary pattern $(0,0)$ the five coarse rows are
$0,1,2,5,8$ and their labels are $4,1,3,0,2$.
All $4\binom{9}{2}=144$ groups pass the same check.

\subsection{Order 54}
Use $(p,d,e,w)=(3,3,5,11)$ and the complete data in
Tables~\ref{tab:affine54-hc} and~\ref{tab:affine54-m}.
The binary array has 12 rows. Its eleven columns are obtained from a
normalized Paley Hadamard matrix of order 12. The maps $A_{i,a}$
depend on both the coordinate $i$ and the binary symbol $a$.
Enumerating $12\cdot3^5=2916$ rows by
\eqref{eq:affine-oa} gives nine MOLS of order 54.
There are $4\binom{11}{2}=220$ small affine partitions to check.

\subsection{Order 98}
Use $(p,d,e,w)=(7,2,3,9)$.
The input matrices have the particularly simple form
\[
 A_{i,a}=\begin{bmatrix}W_{i,a,0}&1&0\\
                         W_{i,a,1}&0&1\end{bmatrix}.
\]
Consequently, for $x=(u,v,w)$, the fine output is
$(v,w)+uW_{i,a}+C_{t,i}$.
Tables~\ref{tab:affine98-hc} and~\ref{tab:affine98-w} give all of $H,W,C$.
There are $28\cdot7^3=9604$ expanded rows and therefore seven MOLS.
The general check has a useful two-dimensional interpretation here.
For fixed $i,j,a,b$, set $\delta=W_{i,a}-W_{j,b}\ne0$.
The seven sets
\[
 C_{t,i}-C_{t,j}+\{u\delta:u\in\F_7\},
       \qquad t\in T_{ij}^{ab},
\]
are parallel affine lines partitioning $\F_7^2$.
A requested pair of fine outputs determines its difference, then $t,u$,
then $v,w$. There are $4\binom{9}{2}=144$ such line partitions.

\section{Small orbit certificates for orders 22 and 74}
Generating DM or QDM seeds under coordinate permutations and entry
automorphisms is a standard direct-construction device; examples include
Abel's constructions at orders 34 and 44~\cite[Lemma 4.1]{abel2008}.
The two actions and all seed rows used here are specified below.
\subsection{Order 22 over the cyclic group}\label{sec:cyclic22}
Take $G=\Z_{21}$, $h=1$, and $w=6$, writing $\infty=\infty_0$.
All finite arithmetic in this subsection is modulo 21.
Define the operation
\[
 T(a,b,c,d,e,f)=(16c,16a,16b,16f,16d,16e),\qquad T^3=\mathrm{id},
\]
with multiplication fixing infinity. This is the multiplier-and-coordinate
permutation of the form in~\cite[Theorem 5.4]{abel2008}, here with multiplier 16;
that source traces the method to Todorov's and Abel--Todorov's
small-order constructions~\cite{todorov1989,abeltodorov1993}.
Use the six base rows and five
supplemental rows in the following two matrices:
\[
R=\begin{pmatrix}
\infty&0&0&0&1&4\\
0&1&4&\infty&14&14\\
0&6&16&3&11&8\\
0&4&6&16&12&3\\
0&7&11&14&3&17\\
0&10&3&0&7&11
\end{pmatrix},\qquad
F=\begin{pmatrix}
0&11&19&7&18&5\\
0&18&12&6&9&15\\
0&3&9&8&5&20\\
0&2&13&17&1&18\\
0&5&1&18&20&10
\end{pmatrix}.
\]
For each row $r$ of $R$, append $r,T(r),T^2(r)$, then append the rows
of $F$. This is a $23\times6$ QDM. Each infinity occurs once per column;
all 15 finite difference lists are permutations of $\Z_{21}$.
Both properties follow by enumerating the displayed rows.
Developing by all 21 additive translations and adjoining the
all-infinity row gives $23\cdot21+1=484$ OA rows, hence four MOLS.
All four completed squares are printed together in Appendix~\ref{app:order22}.

\subsection{Order 74 using quaternion coordinate permutations}
Here $G=\Z_{67}$, $h=7$, $w=8$. Index the coordinates, in order, by
\[
 (1,i,j,k,-1,-i,-j,-k),\qquad i^2=j^2=k^2=ijk=-1.
\]
For computation use $ij=k$, $jk=i$, $ki=j$; reversing the order of
two distinct members of $\{i,j,k\}$ changes the sign, and signs multiply.
These are labels in the quaternion group $Q_8$; the entries themselves
are residues modulo 67 or infinity symbols. Let
\[
B=\begin{pmatrix}
\infty_0&1&29&37&0&35&10&22\\
\infty_1&29&37&1&0&50&43&41\\
\infty_2&39&59&36&0&19&15&33\\
\infty_3&59&36&39&0&52&34&48\\
\infty_4&10&22&35&0&40&21&6\\
\infty_5&43&41&50&0&7&2&58\\
\infty_6&25&55&54&0&3&20&44\\
0&33&37&60&5&21&23&43\\
0&65&19&1&11&41&6&64\\
0&29&9&15&51&47&50&40
\end{pmatrix}.
\]
For every row $b$ and every $g\in Q_8$, append the row whose coordinate
$u$ is $b_{gu}$. Finally append the all-zero row. The resulting
$81\times8$ matrix has every infinity once in each column, and each
of its 28 finite column-pair difference lists is a permutation of $\Z_{67}$.
Develop by the 67 additive translations. Fill the hole with the 49 rows
\[
 (\infty_u,\infty_v,\infty_{u+v},\infty_{u+2v},\ldots,
   \infty_{u+6v}),\qquad u,v\in\F_7,
\]
with subscripts computed modulo 7. The $5427+49=5476$ rows give six MOLS.

\section{Order 85: seven squares from a cyclic difference matrix}
\label{sec:cyclic85}

Work in the additive group $\Z_{85}$; all arithmetic in this section is
ordinary arithmetic modulo 85. Cyclic shifts combined with multiplier
orbits also occur in earlier QDM constructions~\cite[Section 3]{abel2008}.
Here the following three cycles give a complete DM.
Put $\beta=47$, so that
$\beta^2=-1$ and $\beta^4=1$, and use these three seven-entry cycles:
\[
\begin{split}
 A_0&=(6,7,33,68,53,5,16),\\
 A_1&=(36,40,48,76,67,43,23),\\
 A_2&=(38,44,56,73,30,63,60).
\end{split}
\]
Cycle indices are reduced modulo seven. Start an $85\times8$ matrix
$D$ with one all-zero row. For each $a=0,1,2$, then $p=0,1,2,3$,
then $t=0,\ldots,6$, append
\begin{equation}\label{eq:cyclic85-row}
 \bigl(0,\ \beta^p A_a[t],\ \beta^p A_a[t+1],\ \ldots,
             \beta^p A_a[t+6]\bigr).
\end{equation}
This order assigns the appended row number $1+28a+7p+t$.
There are $1+3\cdot4\cdot7=85$ rows.

The following small check proves the difference-matrix condition.
Multiplication by the four invertible residues
\[
 H=\{1,47,84,38\}
\]
partitions the 84 nonzero residues into 21 sets $Hx=\{hx:h\in H\}$,
called \emph{multiplicative orbits}. Each orbit has four elements:
the differences $47-1$, $84-1$, and $38-1$ are all coprime to 85,
so none of these three multipliers fixes a nonzero residue.
Let $\rho(x)$ be the least integer in $Hx$, with residues represented
by $0,\ldots,84$. The set of nonzero orbit representatives is
\[
 R=\{1,2,3,4,5,6,7,8,10,12,13,14,15,17,19,21,23,26,28,30,37\}.
\]
Direct evaluation of the printed cycles gives the following four
finite checks, where every member of $R$ appears exactly once:
\begin{align*}
 \{\rho(A_a[t]):0\le a<3,\ 0\le t<7\}&=R,
       \\
 \{\rho(A_a[t+d]-A_a[t]):0\le a<3,\ 0\le t<7\}&=R,
       \qquad d=1,2,3.
\end{align*}
Thus the check requires only forming four lists of 21 integers and
sorting each to compare it with the displayed $R$.

For a pair consisting of column zero and any other column, cycling $t$
runs through all 21 seed entries. The first check and the four choices
of $p$ supply every nonzero residue once. For two nonzero columns,
their cyclic separation is some $d\in\{1,\ldots,6\}$, and their
differences are the four multiples of $A_a[t+d]-A_a[t]$.
The three difference checks cover separations 1, 2, and 3.
The other separations reverse one of these differences after a shift
of $t$; negation does not change its orbit because $-1\in H$.
Consequently every column pair has every nonzero difference once.
The initial all-zero row supplies zero once, proving that $D$ is a
difference matrix.

Its first column is already zero, so the common development rule
\eqref{eq:dm-squares} gives the seven squares directly:
\[
 L_j(r,c)=D_{r,j}+c\pmod{85},\qquad
 1\le j\le7,\quad 0\le r,c<85.
\]
This proves $N(85)\ge7$.

\section{Relative difference matrices for orders 78 and 92}
\label{sec:relative-difference}

The two constructions share the same development and hole filling.
Only their fourteen seed rows differ. They give seven MOLS at each order.

\subsection{Relative development and hole filling}
Let $p=11$ for order 78 and $p=13$ for order 92. Write
$G=\F_7\times\F_p$ additively and $K=\F_7\times\{0\}$.
An index-one \emph{relative difference matrix} over $(G,K)$ is a
matrix whose difference list in any two columns contains each element
of $G\setminus K$ exactly once and no element of $K$.
Thus translating its rows covers pairs from different cosets of $K$,
leaving one seven-symbol hole per coset. This is the multiple-hole
counterpart of the QDM development in Section~\ref{sec:qdm}.

Let $H$ and $Y$ be $14\times9$ arrays over $\F_7$ and $\F_p$.
Call these the \emph{coarse} and \emph{fine} arrays, respectively.
Let $Q_p$ be the nonzero quadratic residues in $\F_p$, and let
$\chi_p(x)$ be $0$ if $x=0$, $+1$ if $x\in Q_p$, and $-1$ otherwise.
The complete finite condition on the seeds is this:
for every $i<j$ and $a\in\F_7$, exactly two rows $r,s$ satisfy
$H_{r,j}-H_{r,i}=H_{s,j}-H_{s,i}=a$, and
\begin{equation}\label{eq:relative-character}
 \chi_p(Y_{r,j}-Y_{r,i})\,
 \chi_p(Y_{s,j}-Y_{s,i})=-1.
\end{equation}
In particular all fine row entries are distinct. There are
$\binom92\cdot7=252$ such tests.

For each seed row $r$ and each $m\in Q_p$, output
\begin{equation}\label{eq:relative-seed}
 D_{r,m,i}=(H_{r,i},mY_{r,i}),\qquad 0\le i<9.
\end{equation}
The two fine differences in each coarse fiber have opposite characters.
Their products by $Q_p$ therefore list the two disjoint cosets of
$Q_p$ in $\F_p^*$, each once. Consequently $D$ is a relative
difference matrix with $14(p-1)/2=7(p-1)$ rows.

Translate every row of $D$ by every element of $G$. The resulting
$49p(p-1)$ rows contain every ordered pair of finite symbols from
different holes exactly once in any two coordinates: their difference
chooses the unique base row, then either prescribed symbol chooses
the translation. The holes are $B_b=\F_7\times\{b\}$, for $b\in\F_p$.
Adjoin one common symbol $\infty$. For each hole use
$E=\F_2[z]/(z^3+z+1)$ and the complete nine-coordinate field array
\begin{equation}\label{eq:relative-fill}
 (x,y,x+a_1y,\ldots,x+a_7y),\qquad x,y\in E,
\end{equation}
where $a_1,\ldots,a_7$ are the seven nonzero field elements in increasing
bit-label order. Encode $e_0+e_1z+e_2z^2$ by $e_0+2e_1+4e_2$.
Map label zero to $\infty$ and label $e\ne0$ to $(e-1,b)$.
Omit the all-zero input $(x,y)=(0,0)$ from each fill, and append
one global all-$\infty$ row.

Any two field coordinate forms in~\eqref{eq:relative-fill} are
independent. Hence the omitted row is the only fill row with infinity
in both of two specified coordinates. All finite pairs in one hole
and all pairs with exactly one infinity occur once in that hole's fill;
the pair $(\infty,\infty)$ occurs once in the final row.
Together with the developed relative matrix, this proves an OA with
\[
 49p(p-1)+63p+1=(7p+1)^2
\]
rows and nine columns, yielding seven MOLS by~\eqref{eq:extract}.
The development of a relative difference matrix into an orthogonal
array with holes is described in \cite[Lemma 2.8]{bailey2020}.
Abel and Li~\cite[Lemma 3.4]{abelli2015} use a relative difference
matrix, followed by adjoining a common point and filling the holes,
to construct three mutually orthogonal diagonal Latin squares of
order~21.

\begin{center}
\begin{tabular}{rrrrr}
\toprule
$n$ & $p$ & relative rows & translated rows & completed OA rows\\
\midrule
78 & 11 & 70 & 5,390 & 6,084\\
92 & 13 & 84 & 7,644 & 8,464\\
\bottomrule
\end{tabular}
\end{center}
Use $(a,b)\mapsto a+7b$ and $\infty\mapsto7p$ for square symbols.
Use increasing seed-row and multiplier order, translations ordered first
by $a$ then $b$, and fills ordered by $b,x,y$; the OA row order itself
does not affect the squares extracted from coordinate addresses.

\subsection{Order 78: quadratic coarse rows and two short sequences}
Index the rows by $(q,a)\in\{0,1\}\times\F_7$ and the columns by
$A,B,T_0,\ldots,T_6$, in that order. Set
\[
 (c_0,c_1)=(1,3),\qquad (b_0,b_1)=(1,2),
\]
and use the sequences over $\F_{11}$
\[
 U_0=(10,8,5,7,4,2,3),\qquad
 U_1=(1,5,4,6,9,7,10).
\]
Define all seed entries by
\begin{align}
 H_{q,a,A}&=0,&H_{q,a,B}&=a-(4c_q)^{-1},&
 H_{q,a,T_t}&=c_q(a+t)^2,\label{eq:relative78-coarse}\\
 Y_{q,a,A}&=0,&Y_{q,a,B}&=b_q,&
 Y_{q,a,T_t}&=U_q(a+t).\label{eq:relative78-fine}
\end{align}
All coarse arithmetic and sequence subscripts are modulo 7;
fine arithmetic is modulo 11.

To verify the coarse multiplicity condition, first observe that the
$A,B$ difference is a permutation of $\F_7$ in each type $q$.
For $A,T_t$ the difference is $c_q(a+t)^2$. Since 1 and 3 have
opposite quadratic characters modulo 7, zero has one root in each
type and any nonzero value has two roots in exactly one type.
Writing $s_q=(2c_q)^{-1}$, the $B,T_t$ difference is
\[
 c_q(a+t-s_q)^2+t,
\]
so the same argument applies. Finally the $T_t,T_u$ difference,
$c_q(u-t)(2a+u+t)$, is a permutation in each type when $t\ne u$.
These cases cover all 36 column pairs.

The fine tests~\eqref{eq:relative-character} are direct finite checks
on the two sequences. A compact way to check them is to test, for
each type and $h=1,2,3$, opposite characters within
\[
 \bigl(U_q(h),U_q(-h)\bigr),\qquad
 \bigl(U_q(s_q+h)-b_q,U_q(s_q-h)-b_q\bigr).
\]
Then check the fixed pair $(b_0,b_1)$, the two vertex pairs
\[
 (U_0(0),U_1(0))=(10,1),\qquad
 (U_0(4)-1,U_1(6)-2)=(3,8),
\]
and the following 21 pairs, for $d=1,2,3$ and $a\in\F_7$:
\begin{equation}\label{eq:relative78-moving}
 \bigl(U_0(a+d)-U_0(a),\quad
       U_1(5a+3d)-U_1(5a+2d)\bigr).
\end{equation}
Each has one quadratic residue and one nonresidue modulo 11.
For moving columns, equality of coarse differences pairs the type-zero
row $a$ with the type-one row $5a+2d$. Cyclic shifts reindex $a$;
reversing a column pair negates both differences and preserves
opposition of characters. Thus these checks imply all 252 tests.
Alternatively, check all 252 conditions directly.

\subsection{Order 92: a restricted coarse matrix}
Use the full seed tables $H,Y$ in Appendix~\ref{app:relative-data},
over $\F_7$ and $\F_{13}$ respectively, and
$Q_{13}=\{1,3,4,9,10,12\}$ in~\eqref{eq:relative-seed}.
A generalized Hadamard matrix $\mathrm{GH}(7,2)$ is, in our orientation,
a $14\times14$ array over $\F_7$ whose column differences contain
each residue twice. The coarse $H$ retains columns
$0,1,2,3,4,6,7,11,13$ of the second matrix in Lampio's
classification supplement~\cite{lampio2015}.
The displayed $H$ has two occurrences of each coarse difference;
together with $Y$, it satisfies all 252 character conditions
in~\eqref{eq:relative-character}.

\section{Order 90: eight squares from a cyclic column action}
\label{sec:cyclic90}

We give a ten-column QDM over the additive group $\F_3^4$ with a nine-symbol hole.

Use the fields
\[
 F=\F_3[u]/(u^2+1),\qquad
 K=F[v]/\bigl(v^2-(1+u)\bigr).
\]
All scalar coefficients are reduced modulo three.
Since $(1+u)^4=-1$ in $F$, the element $1+u$ is a nonsquare;
thus these are fields of orders 9 and 81.
The integer label
$a_0+3a_1+9b_0+27b_1$ represents
$a_0+a_1u+b_0v+b_1uv\in K$.
Addition is coordinatewise modulo three; multiplication is polynomial
multiplication using $u^2=-1$ and $v^2=1+u$.
Put
\[
 \beta=31,\qquad \alpha=62=-\beta,\qquad
 W=\langle1,v\rangle_{\F_3}
   =\{0,1,2,9,10,11,18,19,20\}.
\]
Direct arithmetic gives $\operatorname{ord}(\beta)=10$ and
$\operatorname{ord}(\alpha)=5$.
Write $w_0,\ldots,w_8$ for $W$ in the displayed order, and give
$\infty_h$ the integer label $81+h$. The intercept data are
\begin{center}
\begin{tabular}{c|rrrrrrrrrr}
$c$ & 0&1&2&3&4&5&6&7&8&9\\
\hline
$z_c$ &0&54&59&47&51&0&69&13&4&40
\end{tabular}
\end{center}
The dummy value $z_0$ is unused.

Columns and intercept subscripts are taken modulo ten.
Output nine wholly finite rows $B_w$, one for each $w\in W$, and
90 punctured rows $R_{w,t}$, one for each $w=w_h\in W$ and
$t\in\{0,\ldots,9\}$:
\begin{align}
 B_{w,c}&=w\beta^c,\nonumber\\
 R_{w,t,c}&=
 \begin{cases}
   \infty_h,&c=t,\\
   \alpha^t\bigl(w\beta^{c-t}+z_{c-t}\bigr),&c\ne t.
 \end{cases}
 \label{eq:cyclic90-rows}
\end{align}
The resulting matrix has 99 rows and ten columns.
Each punctured row has exactly one infinity, and each infinity occurs
once in each column.

The finite-difference condition can be checked in the additive
quotient $K/W$. Its nine cosets are labelled by
\[
 \pi(a_0+a_1u+b_0v+b_1uv)=a_1+3b_1.
\]
These labels encode two coordinates over $\F_3$; no multiplication on
the quotient is needed. For the pair of columns $0,d$, put
$L_d=1-\beta^d$. The full rows supply precisely the nine differences
in $L_dW$. For $t\notin\{0,d\}$, the nine punctured rows give
\[
 R_{w,t,0}-R_{w,t,d}
 =(-1)^t wL_d+\alpha^t(z_{-t}-z_{d-t}),
\]
because $\alpha^t\beta^{-t}=(-1)^t$.
They therefore supply one complete affine coset of $L_dW$.
After dividing by $L_d$, its coset label is
\begin{equation}
 q_{d,t}=\pi\!\left(
 \frac{\alpha^t(z_{-t}-z_{d-t})}{1-\beta^d}\right).
 \label{eq:cyclic90-quotient}
\end{equation}
For each table row, list the valid $t$ increasingly in
$\{0,\ldots,9\}\setminus\{0,d\}$:
\begin{center}
\begin{tabular}{ccl}
\toprule
$d$ & $L_d$ & Quotient labels $q_{d,t}$\\
\midrule
1 & 60 & $4,6,3,2,1,7,5,8$\\
2 & 53 & $4,7,6,5,1,8,2,3$\\
3 & 48 & $6,3,7,4,8,1,5,2$\\
4 & 59 & $3,2,8,5,6,4,1,7$\\
5 & 2  & $2,7,3,4,1,5,6,8$\\
\bottomrule
\end{tabular}
\end{center}
Each row is a permutation of the eight nonzero cosets. Together with
the full-row zero coset, these partition $K$ and give every finite
difference exactly once.

The five displayed column pairs suffice. Simultaneously replace a column
$c$ by $c+1$, a full row $B_w$ by $B_{-w}$, and a punctured row
$R_{w,t}$ by $R_{w,t+1}$. Equation~\eqref{eq:cyclic90-rows} shows that
every finite entry is multiplied by $\alpha$, with infinities fixed.
Thus column shifts preserve the difference condition. Interchanging
two columns negates their differences. A shift and, if needed, an
interchange reduce any distinct pair to $(0,d)$ with $1\le d\le5$.
All 45 pairs can also be checked directly.

Develop these 99 rows by all 81 additive translations of $K$, and fill
the hole using the 81 rows
\[
 (\infty_x,\infty_y,\infty_{x+y},\infty_{x+2y},\ldots,
   \infty_{x+8y}),\qquad x,y\in F,
\]
where labels $0,\ldots,8$ encode the field $F$ and each expression
$x+ay$ uses field arithmetic. Here the infinity subscript is the
integer label of that field element. The resulting
$99\cdot81+9^2=8100$ rows form an $\OA(8100,10,90,2)$, giving
eight MOLS of order 90 by the first-two-coordinate addressing rule.
The choice of $W$ is essential to this particular certificate: it is
an additive plane, not a subfield, since $v\in W$ but $v^2=1+u\notin W$.
The proof requires only that $W$ be closed under addition and negation.

\section{Order 60: six squares from affine coordinate operations}\label{sec:order60}

We construct a $68\times8$ quasi-difference matrix over
$G=\Z_{26}\times\Z_2$, with eight additional symbols
$\infty_0,\ldots,\infty_7$. Encode $(t,e)\in G$ as $t+26e$,
and encode $\infty_h$ as $52+h$. Thus the finite labels are
$0,\ldots,51$; their addition is componentwise, not addition modulo 52.
Every translation fixes the infinity symbols.

\paragraph{Three coordinate operations.}
For an eight-coordinate row $v$, define operations $T_0,T_1,T_2$ by
\[
 (T_jv)_i=v_{i\xor 2^j}+c_{j,i},\qquad 0\le i<8,
\]
where $\xor$ denotes bitwise XOR on the coordinate indices, and
\[
\begin{aligned}
 c_0&=(0,0,0,0,0,0,0,0),\\
 c_1&=(26,39,0,13,26,39,0,13),\\
 c_2&=(39,13,39,13,0,26,0,26).
\end{aligned}
\]
Here each $c_{j,i}$ is a label in $G$, and the addition fixes any
infinity entry. For $a=a_0+2a_1+4a_2$, with $a_j\in\{0,1\}$, put
\[
 T(a)=T_0^{a_0}T_1^{a_1}T_2^{a_2}.
\]
The rightmost operation acts first: apply $T_2$ when $a_2=1$, then
$T_1$ when $a_1=1$, and finally $T_0$ when $a_0=1$.
\paragraph{The finite certificate.}
Take eight seeds $s_h=(\infty_h,0,b_{h,0},\ldots,b_{h,5})$, with
\[
\begin{array}{c|rrrrrr}
h&b_{h,0}&b_{h,1}&b_{h,2}&b_{h,3}&b_{h,4}&b_{h,5}\\\hline
0&40&45&26&28&38&39\\
1&23&39&29&5&26&32\\
2&13&5&30&13&22&41\\
3&33&43&22&14&33&8\\
4&47&50&21&25&1&3\\
5&16&22&11&33&37&49\\
6&32&34&31&37&2&20\\
7&24&10&20&4&21&25
\end{array}
\]
and the finite seed
\[
 f=(0,0,14,14,16,16,9,9).
\]
The rows of $Q$ are $T(a)s_h$, first in increasing $h=0,\ldots,7$
and then in increasing $a=0,\ldots,7$, followed by $T(a)f$ for
$a=0,2,4,6$ in that order. This gives $8\cdot8+4=68=52+2\cdot8$
rows. Each of the first 64 rows has one infinity symbol, and the last
four have none. Since the coordinate permutations are XOR by
$1,2,4$, the infinity in $T(a)s_h$ lies in coordinate $a$.
Consequently every infinity symbol occurs once in every column.

For each of the $\binom82=28$ pairs of distinct columns, deleting the
16 rows containing an infinity in either column leaves 52 differences.
Direct substitution of the displayed data, using subtraction in
$\Z_{26}\times\Z_2$, gives each of the labels $0,\ldots,51$ exactly
once in every such list. Thus $Q$ satisfies all three QDM conditions
of Section~\ref{sec:qdm}.

\paragraph{Development and the eight-symbol hole.}
Translate each row of $Q$ by every $g\in G$, leaving its infinity
entry fixed. This gives $68\cdot52=3536$ rows. To fill the remaining
hole, use
\[
 F=\F_2[z]/(z^3+z+1),\qquad
 \enc(x_0+x_1z+x_2z^2)=x_0+2x_1+4x_2.
\]
For every $x,y\in F$, append the row whose coordinate $i$, with $i$
interpreted as a field label, is
\[
 \infty_{\enc(y+ix)},\qquad 0\le i<8.
\]
For distinct $i,j$, the pair $(y+ix,y+jx)$ determines $x$ and $y$
uniquely because $j-i\ne0$. These 64 rows therefore form
$\OA(64,8,8,2)$ on the infinity symbols. Enumerate $x$ first and
$y$ second, in increasing label order, to agree with the implementation.
Together the rows form $\OA(3600,8,60,2)$, and
\eqref{eq:extract} gives six MOLS of order 60. 

\section{Order 225: nine squares from two coordinate functions}\label{sec:order225}
Use the additive group $G=\Z_{15}\times\Z_{15}$, with
$\enc(u,v)=u+15v$. All arithmetic on pairs is componentwise modulo~15,
not modulo~225. Define $B,C:\Z_{15}\to\Z_{15}$ by
\begin{center}\small
\setlength{\tabcolsep}{4pt}
\begin{tabular}{r*{15}{r}}
\toprule
$t$ &0&1&2&3&4&5&6&7&8&9&10&11&12&13&14\\
\midrule
$B(t)$ &0&4&6&8&13&3&5&14&1&10&12&2&7&9&11\\
$C(t)$ &0&12&10&13&7&6&4&1&14&11&9&8&2&5&3\\
\bottomrule
\end{tabular}
\end{center}

Form a $225\times10$ matrix $D$, with row $(x,y)\in G$ at position
$r=x+15y$. Write $D_i(x,y)=D_{r,i}$. The two components of each entry
are given below and reduced modulo~15.
\begin{center}
\begin{tabular}{cll}
\toprule
Column $i$ & First component of $D_i(x,y)$ & Second component\\
\midrule
0 & $B(x)+B(y)$ & $x+y$\\
1 & $C(x)+y$ & $x+B(y)$\\
2 & $B(y)$ & $B(x)+B(y)$\\
3 & $y$ & $C(x)+y$\\
4 & $x+C(y)$ & $x+C(y)$\\
5 & $C(y)$ & $x$\\
6 & $B(x)+y$ & $B(x)$\\
7 & $B(x)+C(y)$ & $C(x)+B(y)$\\
8 & $C(x)+C(y)$ & $B(x)+y$\\
9 & $C(x)+B(y)$ & $C(x)$\\
\bottomrule
\end{tabular}
\end{center}
For example, row $(1,2)$ has index~31, with $D_0(1,2)=(10,3)$
and $D_1(1,2)=(14,7)$.

\paragraph{Why the differences are permutations.}
Let $\operatorname{id}(t)=t$. Direct inspection of the two short tables
shows that
\[
 \operatorname{id},\quad B,\quad C,\quad
 \operatorname{id}-B,\quad\operatorname{id}-C,\quad B-C
\]
are permutations of $\Z_{15}$. Thus the difference of any two distinct
functions in $\{0,\operatorname{id},B,C\}$ is a permutation.

Each formula for $D_i$ has four function slots: the $x$- and $y$-terms
in its first and second components, with missing terms read as zero.
Of the 45 pairs of columns, 40 agree in exactly one slot. In their
difference that slot cancels, so, after possibly exchanging the input
variables or the two output components, the difference is
\[
 (u,v)=\bigl(f(y),\,g(x)+h(y)\bigr),
\]
where $f$ and $g$ are differences of distinct functions from the displayed
set. Recover $y=f^{-1}(u)$ and then $x=g^{-1}(v-h(y))$ uniquely.

Of the five remaining pairs, $D_4-D_5=(x,C(y))$ is immediate. The other
four are inverted as follows, with $(u,v)$ the indicated difference.
\begin{center}\small
\begin{tabular}{lll}
\toprule
Difference & Recover one variable from & Then recover the other from\\
\midrule
$D_2-D_4$ & $B(x)=v-u$ & $(B-C)(y)=u+x$\\
$D_3-D_4$ & $C(x)=v-u$ & $(\operatorname{id}-C)(y)=u+x$\\
$D_6-D_4$ & $y=u-v$ & $(B-\operatorname{id})(x)=v+C(y)$\\
$D_9-D_4$ & $B(y)=u-v$ & $(C-\operatorname{id})(x)=v+C(y)$\\
\bottomrule
\end{tabular}
\end{center}
All the functions inverted are among the permutations already checked,
or their negatives. This proves the difference condition for every
column pair. Individual columns need not be permutations: column~4,
for example, has equal components.

\paragraph{Development.}
Subtract $D_0(x,y)$ from every entry in row $(x,y)$ and apply
\eqref{eq:dm-squares}, using $r=x+15y$ and increasing encoded group
labels $g=0,\ldots,224$. The resulting nine squares have order~225,
so $N(225)\ge9$. 

\paragraph{Provenance.}
The tables $B,C$ were obtained by a common relabeling of the inputs of
three maps extracted from the order-15 MOLS in SageMath's database,
whose \code{MOLS_15_4} entry credits Ian Wanless~\cite{sage225}.

\section{Larger affine constructions}\label{sec:finite-affine}
The six orders in Table~\ref{tab:finite-parameters} use the same partition
principle as the smaller affine constructions. Four produce difference
matrices. The remaining two are translated-array instances of the general
affine array construction: each fine output is a column-dependent linear
map plus a common translation. 

\begin{table}[ht]
\centering\small
\begin{tabular}{rllrrrrrrr}
\toprule
$n$ & Route & Coarse group & $p$ & $d$ & $a$ & $\lambda$ & Rows & Columns & MOLS\\
\midrule
144 & DM & $\Z_3^2$ & 2 & 4 & 3 & 2 & 18 & 15 & 14\\
162 & Translated array & $\{0,1\}$ & 3 & 4 & 3 & 3 & 12 & 11 & 9\\
192 & DM & $\Z_3$ & 2 & 6 & 4 & 4 & 12 & 10 & 9\\
288 & DM & $\Z_3^2$ & 2 & 5 & 4 & 2 & 18 & 18 & 17\\
324 & DM & $\Z_2^2$ & 3 & 4 & 3 & 3 & 12 & 12 & 11\\
486 & Translated array & $\{0,1\}$ & 3 & 5 & 4 & 3 & 12 & 11 & 9\\
\bottomrule
\end{tabular}
\caption{Fine vectors lie in $\F_p^d$, inputs in $\F_p^a$, and
$\lambda=p^{d-a}$. ``Rows'' and ``Columns'' refer to the coarse array $H$.}
\label{tab:finite-parameters}
\end{table}

In this section a matrix $M_i$ always has $d$ output rows and $a$ input
columns. Its entries are residues modulo $p$; all vector sums are
coordinatewise modulo $p$. A vector $(v_0,\ldots,v_{d-1})$ is encoded as
$\sum_b v_bp^b$. If the coarse alphabet has size $s$, the full symbol is
encoded as $h+s\enc_p(v)$. Thus the coarse digit is the least significant
part. Tables~\ref{tab:finite-c144}--\ref{tab:finite-c486} give the offsets.

\subsection{Four difference matrices}
For each of orders 144, 192, 288, and 324, form one row for every coarse
row $t$ and every input vector $y\in\F_p^a$, with entry
\begin{equation}\label{eq:finite-dm}
 D_{(t,y),i}=\bigl(H_{t,i},M_i y+c_{t,i}\bigr).
\end{equation}
Use increasing $t$ and then increasing little-endian encoding of $y$ as the
row order. Apply the normalization and development algorithm of
Section~\ref{sec:affine-lifts}. The finite check is exactly its difference
partition condition: for each column pair and each coarse difference, the
$\lambda$ translated images of $M_i-M_j$ partition $\F_p^d$.

\paragraph{Orders 144 and 288: one coarse array.}
Encode $(a,b)\in\F_3^2$ by $a+3b$. Define, with arithmetic modulo three,
\[
 P(a,b)=(a^2+b^2,\;2ab+2b^2),\qquad J(a,b)=(b,\;a+2b).
\]
For $u,v=0,\ldots,8$, let $Q(u,v)=P(\operatorname{digits}_3(u)
-\operatorname{digits}_3(v))$. Build the $18\times18$ array
\begin{equation}\label{eq:finite-h9}
 H=\begin{pmatrix}Q&JQ\\-Q&J^2Q\end{pmatrix}.
\end{equation}
Here $J$, $J^2$, and negation act on each entry. List the top block before
the bottom block, and within a block use increasing $u$; list the left
block of columns before the right block. Retain columns 0 through 14 for
order 144, and all columns for order 288. Each coarse column difference
lists every element of $\F_3^2$ twice.

For a $d$-bit integer $x$, define the linear shift operation
\begin{equation}\label{eq:finite-shift}
 T_{d,f}(x)=((2x)\bmod2^d)\mathbin{\oplus}
 \begin{cases}f,&x\ge2^{d-1},\\0,&x<2^{d-1},\end{cases}
\end{equation}
where $\oplus$ denotes bitwise exclusive-or. The encoded columns of $M_i$ are
\begin{equation}\label{eq:finite-binarymaps}
 \bigl(i,T_{d,f}(i),\ldots,T_{d,f}^{a-1}(i)\bigr).
\end{equation}
Use $(d,a,f)=(4,3,3)$ at order 144 and $(5,4,5)$ at order 288. Matrix-vector
multiplication means taking the XOR of those encoded columns selected by
the 1-bits of $y$. The offsets in Tables~\ref{tab:finite-c144} and
\ref{tab:finite-c288} satisfy the binary index-two criterion of
Section~\ref{sec:affine-lifts}, giving fourteen and seventeen MOLS,
respectively.

\paragraph{Order 192: a codimension-two lift.}
Use the $12\times10$ ternary coarse array in Table~\ref{tab:finite-h192}.
It has index four: every coarse difference occurs four times. The array is
extracted from the published 36-row ternary array in Sloane's
library~\cite{sloane36}. Set $(p,d,a)=(2,6,4)$, and use
\eqref{eq:finite-binarymaps} with $f=3$ and $i=0,\ldots,9$.
The six-bit offsets in Table~\ref{tab:finite-c192} finish the construction.
The output group is $\Z_3\times\Z_2^6$.

\paragraph{Order 324: a ternary lift.}
Use the $12\times12$ array in Table~\ref{tab:finite-h324}, over
$\Z_2^2$. Its symbols are two-bit integers, so coarse subtraction is XOR.
This ingredient is printed by Luo, Zhang, and Chen~\cite[p.~740]{luo2017},
who attribute it to earlier work of Zhang. Each coarse difference occurs
three times. Set $p=3$, $d=4$, $a=3$ and define
\[
 T(a,b,c,d)=(d,\;a+2d,\;b,\;c)\pmod3.
\]
Let $\mu_i$ be the four ternary digits of $i$, for $0\le i<12$.
The three columns of $M_i$ are $\mu_i,T\mu_i,T^2\mu_i$.
Use the offsets of Table~\ref{tab:finite-c324}. Equation~\eqref{eq:finite-dm}
then gives the difference matrix over $\Z_2^2\times\Z_3^4$.

\subsection{Orders 162 and 486: translated arrays}
The common binary coarse array $B$ is in Table~\ref{tab:finite-binary}.
Every ordered pair of bits occurs three times in every pair of its eleven
columns. It is a normalized Hadamard-derived array of Paley
type~\cite{paley1933}.

At order 486 take $H=B$, $p=3$, $d=5$, $a=4$. Define
\[
 \begin{aligned}
 T(v_0,v_1,v_2,v_3,v_4)&=(2v_4,\;v_0+v_4,\;v_1,\;v_2,\;v_3)\pmod3,\\
 M_i y&=T^i(y_0,y_1,y_2,y_3,0),\qquad 0\le i<11.
 \end{aligned}
\]
Use Table~\ref{tab:finite-c486} for the offsets.
At order 162 take $p=3$, $d=4$, $a=3$, and reorder the columns of $B$ as
\[
 (2,5,6,4,10,8,3,9,7,0,1).
\]
The maps and offsets are Tables~\ref{tab:finite-m162} and
\ref{tab:finite-c162}.

For either order, iterate through every $t$, every $y\in\F_3^a$, and
every $z\in\F_3^d$, in that order and with vector encodings increasing.
Output the eleven entries
\begin{equation}\label{eq:finite-translated}
 A_i(t,y,z)=H_{t,i}+2\enc_3\bigl(z+M_i y+c_{t,i}\bigr).
\end{equation}
The pair $(A_0,A_1)$ identifies a grid cell, and $A_2,\ldots,A_{10}$
give its entries in the nine squares. This is the general affine OA
construction of Section~\ref{sec:affine-arrays}, with fine input $(y,z)$
and map $[M_i\mid I_d]$; the maps do not depend on the coarse bit here.

Its finite condition has a particularly short form. For each pair of
columns $i,j$ and each ordered coarse pair $(h,h')$, take the three rows
$t$ with $(H_{t,i},H_{t,j})=(h,h')$. The sets
\[
 \{(M_i-M_j)y+c_{t,i}-c_{t,j}:y\in\F_3^a\}
\]
must partition $\F_3^d$. Given two prescribed fine outputs, their difference
therefore determines $(t,y)$ uniquely, and then either output determines
$z$. This proves the OA pair property and the correctness of the expansion.

An index-one difference matrix over $\Z_2\times\Z_3^d$ with three or more
columns cannot serve here. After normalization, two nonzero columns would
be permutations of the group. Each binary projection would sum to 1 modulo
2, since it contains $3^d$ ones. Their difference would sum to 0, yet would
also have to be a group permutation and sum to 1. The translated-array
route avoids that obstruction.

\section{Assembling local squares on selected plane points}
\label{sec:geometric}
For each of the twelve orders
$175,183,207,328,335,427,430,442,456,469,482,494$,
we select a set of points in a finite plane. Intersecting the plane's
lines with this set gives blocks in which every pair of points occurs
exactly once. Latin squares on these blocks then combine into Latin
squares on the selected point set.

\subsection{Plane incidence and local filling}
The construction separates the placement of entries from the choice of their
values. A finite plane supplies small sets of points, called blocks, so that
each pair of distinct points belongs to exactly one block. On each block we
construct small Latin squares over a finite field. We then use these local
squares to fill the corresponding cells of the global squares.

There are three different sizes to keep track of: the plane has order $q$,
the selected point set $S$ has size $n$, and an individual block has size $s$.
The desired result is $k$ MOLS of order $n$. Plane coordinates use $\F_q$;
the local squares on a block use $\F_s$. These fields need not have the same
characteristic. We first describe the plane and its blocks, then give the
local filling rule and explain why the resulting squares are Latin and
mutually orthogonal.

\paragraph{Points, lines, and directions.}
A \emph{projective plane of order $q$} consists of $q^2+q+1$ points and
equally many subsets called lines. Each line contains $q+1$ points, every
two distinct points lie on exactly one line, and every two distinct lines
meet in exactly one point.

For the ordinary plane, start with the affine grid $\F_q^2$. Its lines are
$y=mx+b$ and the vertical lines $x=a$. Lines with a fixed slope $m$ are
parallel in this grid. Adjoin one point $I_m$ to all lines of that slope,
and one point $I_\infty$ to all vertical lines. Finally, add a line
containing these $q+1$ new points. Parallel lines now meet at their common
infinity point. These infinity points are included in the point set and
may become row labels, column labels, and symbols of the final squares.
Two constructions below change the direction functions $mx$, while
retaining this way of adjoining infinity points.

\paragraph{The field labels used in this section.}
For $q=p^d$, the integer $a_0+pa_1+\cdots+p^{d-1}a_{d-1}$ denotes
$a_0+a_1\theta+\cdots+a_{d-1}\theta^{d-1}$, with the first coefficient
least significant. Add coefficients modulo $p$ and reduce polynomial products
by the following monic polynomial. The same convention applies to local
block fillings.
\begin{center}
\begin{tabular}{rcl}
\toprule Field size & Characteristic & Defining polynomial\\\midrule
4&2&$X^2+X+1$\\
8&2&$X^3+X+1$\\
9&3&$X^2+1$\\
16&2&$X^4+X+1$\\
25&5&$X^2+2$\\
27&3&$X^3+2X+2$\\
32&2&$X^5+X^2+1$\\
\bottomrule
\end{tabular}
\end{center}
All other field sizes below are prime, and use arithmetic modulo that prime.
In an odd field put $\chi(0)=0$; for nonzero $a$ put $\chi(a)=1$ if
$a=b^2$ for some nonzero $b$, and $\chi(a)=-1$ otherwise.
Compute $\chi$ by first listing the nonzero squares.
For $\F_{25}$ write $\alpha=5$ (so $\alpha^2=3$) and define
\[
 \operatorname{Nm}(a+5b)=a^2+2b^2\pmod5.
\]
Thus $\operatorname{Nm}$ can be computed from two base-five digits.
For $\F_{27}$ write $\beta=3$; then $\beta^3=\beta+1$.

\paragraph{Enumerate all points and lines.}
Give affine point $(x,y)$ the integer ID $qx+y$, point $I_m$ the ID
$q^2+m$ for each $m\in\F_q$, and point $I_\infty$ the ID $q^2+q$.
There are $q^2$ affine points and $q+1$ infinity points.
Choose the following direction functions:
\begin{equation}\label{eq:geo-directions}
 f_m(x)=\begin{cases}
 mx,&\text{ordinary plane (all recipes except 456 and 482)},\\
 mx+m^3x^9,&q=27\text{ and order }456,\\
 mx,&q=27,\text{ order }482,\ m=0\text{ or }\chi(m)=1,\\
 mx^3,&q=27,\text{ order }482,\ \chi(m)=-1.
 \end{cases}
\end{equation}
Enumerate the lines in this order:
\begin{enumerate}
\item $\{(x,f_m(x)+b):x\in\F_q\}\cup\{I_m\}$, for $m=0,\ldots,q-1$
      and then $b=0,\ldots,q-1$;
\item $\{(a,y):y\in\F_q\}\cup\{I_\infty\}$, for $a=0,\ldots,q-1$;
\item $\{I_m:m\in\F_q\}\cup\{I_\infty\}$.
\end{enumerate}
This gives $q^2+q+1$ lines, each containing $q+1$ points. The stated order
fixes the enumeration; correctness does not depend on it.

To check unique pair coverage, first consider two affine points $(x,y)$ and
$(x',y')$ with $x\ne x'$. They lie on a line of direction $m$ precisely when
\[
 y-y'=f_m(x)-f_m(x').
\]
If each difference $f_m-f_j$, $m\ne j$, is a permutation of $\F_q$, then
the right-hand sides are distinct as $m$ varies: equality for $m,j$ would
give $(f_m-f_j)(x)=(f_m-f_j)(x')$, contradicting injectivity.
The $q$ right-hand sides therefore run through the entire field, so they
determine exactly one $m$, followed by $b=y-f_m(x)$. Points with the same
first coordinate use their unique vertical line. An affine point and $I_m$
determine the intercept of a line in direction $m$; an affine point and
$I_\infty$ determine a vertical line. Two infinity points lie on the
infinity line.

For the ordinary plane, $f_m-f_j=(m-j)x$ is a permutation whenever
$m\ne j$. For the other two planes, all direction functions are additive:
in characteristic three, $(x+x')^3=x^3+x'^3$, and the same identity holds
for ninth powers. Their differences are therefore linear maps over $\F_3$.
For such a map on a finite field, it suffices to show that only zero maps
to zero; the map is then injective and hence a permutation.

For order 456, write $a=m-j\ne0$. The difference is $ax+a^3x^9$.
A nonzero solution of $ax+a^3x^9=0$ would imply
$-1=a^2x^8$, a square. But $-1$ is a nonsquare in $\F_{27}$:
a square root would have multiplicative order four, whereas
$|\F_{27}^{*}|=26$. For order 482, differences within either exponent
class are nonzero multiples of $x$ or $x^3$, hence permutations.
This also handles a difference involving $m=0$. Between the two nonzero
classes, a nonzero kernel element would give $a/b=x^2$, where one of
$a,b$ is a square and the other is a nonsquare, which is impossible.
Thus all three line enumerations have the required unique-pair property.

The order-456 multiplication is a \emph{generalized twisted field}:
it is Albert's $mx-cm^{\sigma}x^{\tau}$ with $c=-1$,
$\sigma:m\mapsto m^3$, and $\tau:x\mapsto x^9$~\cite[Eq.~(3)]{albert1961}.
The order-482 plane uses Andr\'e's norm-class switching
construction~\cite{andre1954,andre}: the norm $m^{13}$ from $\F_{27}$ to
$\F_3$ distinguishes the square and nonsquare classes, and the latter
uses the field automorphism $x\mapsto x^3$.

\paragraph{Restrict the lines to a selected point set.}
Each recipe below specifies a set $S$ of $n$ plane points. Sort their IDs
and relabel them $0,\ldots,n-1$; these will be the row labels, column labels,
and symbols of the final squares. Intersect each line with $S$, convert
its points to the new labels, sort them, and retain the intersection if
it has at least two points. Call each retained subset a \emph{block}.
Its size $s$ may be smaller than $q+1$ and may differ between blocks.

Any two selected points were together on exactly one plane line, so they
are together in exactly one retained block. A collection of blocks with
this property is called a \emph{pairwise balanced design} (PBD).
In particular, distinct blocks share at most one point. To fill a global
cell $(x,y)$ with $x\ne y$, we can therefore use the unique block containing
$x$ and $y$. Blocks sharing a point $x$ will share the diagonal cell
$(x,x)$, which is why the local squares must agree there.
Record the sizes of all line intersections, including zero and one;
these give the checks in Table~\ref{tab:geo-histograms}.

Selecting plane points and filling the resulting PBDs is a classical route
to small-order MOLS. Brouwer's separable designs gave, for example, eight
MOLS of order 189, twelve of order 253, and nine of order
357~\cite{brouwer1980}. Greig's projective-plane configurations gave
seven MOLS of order 158 and nine of order 254~\cite{greig1999}.
The following recipes choose point sets whose block sizes support the
required local squares.

\paragraph{Construct the local squares.}
Let $B=(b_0,\ldots,b_{s-1})$ be a block in increasing order, where $s$ is
a prime power. Identify local index $u$ with the element of $\F_s$ having
integer label $u$. The correspondence $u\mapsto b_u$ converts a local
field symbol into a global point label; arithmetic is performed on the
local symbols, not on the point labels.

Suppose $s\ge k+2$. Use the $k$ field elements with integer labels
$2,\ldots,k+1$ as coefficients, and define
\begin{equation}\label{eq:geo-idempotent}
 Q_j(u,v)=(1-t)u+tv,\qquad t=j+1,\quad j=1,\ldots,k.
\end{equation}
The coefficient $t$ is a field label, not repeated addition of the field
unit. These choices are distinct and avoid both zero and one. For fixed
$u$, varying $v$ gives every field element because $t\ne0$; for fixed $v$,
varying $u$ does the same because $1-t\ne0$. Thus each $Q_j$ is Latin.
Also $Q_j(u,u)=u$. This property is called \emph{idempotence}; after
relabeling, it gives the value $b_u$ at diagonal cell $(b_u,b_u)$ in
every block containing that point.

For orthogonality, take distinct coefficients $t,t'$ and prescribe their
two outputs $a,b\in\F_s$. The equations are
$a=u+t(v-u)$ and $b=u+t'(v-u)$. They have the unique solution
\[
 d=\frac{b-a}{t'-t},\qquad u=a-td,\qquad v=u+d.
\]
The denominator is nonzero, so every ordered output pair occurs exactly
once. In particular, equal outputs require $d=0$ and hence $u=v=a$;
the two squares agree only on their common diagonal.

\paragraph{Assemble the global squares.}
Initialize $L_j(x,x)=x$ for every global point $x$ and every
$j=1,\ldots,k$. For each block $B$ and each ordered pair of distinct
local indices $u,v$, assign
\begin{equation}\label{eq:geo-assembly}
 L_j(b_u,b_v)=b_{Q_j(u,v)}\quad(j=1,\ldots,k).
\end{equation}
Every off-diagonal cell is assigned exactly once, because its row and
column labels determine a unique block. The diagonal has already been
filled with the values required by every local square.

Fix a global row $x$. The sets $B\setminus\{x\}$, over blocks containing
$x$, partition all other points. In the local row for $x$, idempotence
puts $x$ on the diagonal, and the remaining entries give each symbol of
$B\setminus\{x\}$ once. Combining these rows therefore gives every global
symbol once. The same argument applies to columns.

Now consider two global squares and prescribe unequal outputs $a,b$.
Their unique common block is the only block that can produce both
symbols. Within that block, local orthogonality determines exactly one
cell, and the cell is off the diagonal because its two outputs differ.
For equal outputs $(a,a)$, the local argument above excludes every
off-diagonal cell, while the global cell $(a,a)$ supplies the pair.
Thus the assembled squares are mutually orthogonal. This is the
classical idempotent PBD filling construction of Bose, Shrikhande, and
Parker~\cite{boseshrikhandeparker1960}; see also the MOLS survey of
Colbourn and Dinitz~\cite[Section 4]{colbourndinitz2001}.

\paragraph{A worked example: two MOLS of order 13.}
Take the ordinary projective plane over $\F_3$ and select all its points.
Thus $q=3$, the final order is $n=3^2+3+1=13$, and every block has
size $s=4$. The nine affine points have labels $3x+y$ for
$x,y\in\{0,1,2\}$; the points $I_0,I_1,I_2,I_\infty$ have labels
$9,10,11,12$, respectively. The restriction step leaves
every line unchanged, and no relabelling is needed. The thirteen lines,
written as sorted blocks in the stated enumeration order, are
% BEGIN_PLANE_EXAMPLE_BLOCKS
\[
\begin{aligned}
&(0,3,6,9),\quad (1,4,7,9),\quad (2,5,8,9)\\
&(0,4,8,10),\quad (1,5,6,10),\quad (2,3,7,10)\\
&(0,5,7,11),\quad (1,3,8,11),\quad (2,4,6,11)\\
&(0,1,2,12),\quad (3,4,5,12),\quad (6,7,8,12)\\
&(9,10,11,12)
\end{aligned}
\]
% END_PLANE_EXAMPLE_BLOCKS
For example, the line $y=2x+1$ gives $(1,3,8,11)$: its three affine
points are $(0,1),(1,0),(2,2)$, and its fourth point is $I_2$.
Every unordered pair of labels occurs in exactly one of these blocks.

The local fillings use $\F_4$, even though the plane was defined over
$\F_3$. Write $\F_4=\F_2[\theta]/(\theta^2+\theta+1)$, with labels
$0,1,2=\theta,3=1+\theta$. Addition is two-bit XOR; the nontrivial
products are $2\cdot2=3$, $2\cdot3=1$, and $3\cdot3=2$.
For $k=2$, the two slopes are $t=2,3$, giving
$Q_1(u,v)=3u+2v$ and $Q_2(u,v)=2u+3v$. Their complete tables are
% BEGIN_PLANE_EXAMPLE_LOCAL
\[
Q_1=\begin{pmatrix}
0&2&3&1\\
3&1&0&2\\
1&3&2&0\\
2&0&1&3
\end{pmatrix},\qquad Q_2=\begin{pmatrix}
0&3&1&2\\
2&1&3&0\\
3&0&2&1\\
1&2&0&3
\end{pmatrix}.
\]
% END_PLANE_EXAMPLE_LOCAL
Rows, columns, and entries of these tables are local field labels
$0,1,2,3$. To use a table on a block $(b_0,b_1,b_2,b_3)$, replace each
local label $i$ by the corresponding global point label $b_i$.

For instance, global cell $(1,11)$ belongs to block $(1,3,8,11)$,
so its local coordinates are $(u,v)=(0,3)$. The table entries are
$Q_1(0,3)=1$ and $Q_2(0,3)=2$, hence
\[
 L_1(1,11)=b_1=3,\qquad L_2(1,11)=b_2=8.
\]
Conversely, the output pair $(3,8)$ identifies this same block, where
its local output pair is $(1,2)$. Adding the two local equations gives
$u+v=3$; substitution in $3u+2v=1$ gives $u=0$, then $v=3$.
The pair therefore occurs only at global cell $(1,11)$.
The infinity points are treated in exactly the same way: in block
$(9,10,11,12)$, cell $(10,11)$ has local coordinates $(1,2)$, so
its two outputs are $9$ and $12$.

Set $L_j(x,x)=x$ and apply the same rule to every off-diagonal cell
of every block. This gives the following two squares, with rows and
columns in order $0,\ldots,12$:
% BEGIN_PLANE_EXAMPLE_GLOBAL
\begingroup
\small
\setcounter{MaxMatrixCols}{13}
\setlength{\arraycolsep}{1.7pt}
\[
L_1=\begin{pmatrix}
0&2&12&6&8&7&9&11&10&3&4&5&1\\
12&1&0&8&7&6&10&9&11&4&5&3&2\\
1&12&2&7&6&8&11&10&9&5&3&4&0\\
9&11&10&3&5&12&0&2&1&6&7&8&4\\
10&9&11&12&4&3&2&1&0&7&8&6&5\\
11&10&9&4&12&5&1&0&2&8&6&7&3\\
3&5&4&9&11&10&6&8&12&0&1&2&7\\
5&4&3&10&9&11&12&7&6&1&2&0&8\\
4&3&5&11&10&9&7&12&8&2&0&1&6\\
6&7&8&0&1&2&3&4&5&9&11&12&10\\
8&6&7&2&0&1&5&3&4&12&10&9&11\\
7&8&6&1&2&0&4&5&3&10&12&11&9\\
2&0&1&5&3&4&8&6&7&11&9&10&12
\end{pmatrix},\qquad L_2=\begin{pmatrix}
0&12&1&9&10&11&3&5&4&6&8&7&2\\
2&1&12&11&9&10&5&4&3&7&6&8&0\\
12&0&2&10&11&9&4&3&5&8&7&6&1\\
6&8&7&3&12&4&9&10&11&0&2&1&5\\
8&7&6&5&4&12&11&9&10&1&0&2&3\\
7&6&8&12&3&5&10&11&9&2&1&0&4\\
9&10&11&0&2&1&6&12&7&3&5&4&8\\
11&9&10&2&1&0&8&7&12&4&3&5&6\\
10&11&9&1&0&2&12&6&8&5&4&3&7\\
3&4&5&6&7&8&0&1&2&9&12&10&11\\
4&5&3&7&8&6&1&2&0&11&10&12&9\\
5&3&4&8&6&7&2&0&1&12&9&11&10\\
1&2&0&4&5&3&7&8&6&10&11&9&12
\end{pmatrix}.
\]
\endgroup
% END_PLANE_EXAMPLE_GLOBAL
Each row and column contains all thirteen symbols, and the superposition
contains every ordered pair once. In particular, equal pairs $(x,x)$
occur at the diagonal, while the block construction supplies the other
$13\cdot12=156$ pairs.

\paragraph{The common-center exception.}
Orders 427 and 482 also have blocks of size $s=k+1$. Such a block has only
$s-2=k-1$ coefficients available in the idempotent formula, one fewer
than needed. A different field formula supplies $s-1=k$ local squares,
but changes the diagonal. We can accommodate that change because all
these smaller blocks contain one point $c$, their \emph{center}, and
any two of them meet only at $c$. Such a collection is called a
\emph{brush}. Each noncentral point then belongs to at most one of the
smaller blocks.

The corresponding filling method is due to Wallis and
Zhu~\cite{walliszhu1983}; Zhu used it to construct six MOLS of order
69~\cite[Theorems 1--3]{zhu1984}. It is also a special case of
near-clear-set filling~\cite[Section 4]{colbourndinitz2001}.
Use the following modification:
\begin{enumerate}
\item Initialize the global diagonal and fill all ordinary off-diagonal
      block cells using~\eqref{eq:geo-idempotent}--\eqref{eq:geo-assembly}.
\item On each exceptional block put $c$ first, giving it local label zero,
      and put the other labels in increasing order. Use
      $Q_j(u,v)=u+tv$ with $t=j$, $j=1,\ldots,k$, and apply
      the assignment~\eqref{eq:geo-assembly} to \emph{all} local cells,
      including the diagonal.
\end{enumerate}
The coefficients now range over all nonzero elements of $\F_s$.
Each local square is Latin, and two distinct coefficients recover a
prescribed output pair by
$v=(b-a)/(t'-t)$ and $u=a-tv$. They are therefore mutually orthogonal.
All exceptional squares give zero at local cell $(0,0)$, so their
assignments agree at the one shared global cell $(c,c)$.

For a noncentral point in an exceptional block $B$, that block now fills
its entire local row and column, supplying exactly the symbols of $B$.
All other incident blocks supply the disjoint sets outside $B$; their
diagonal cells were not assigned by the ordinary filling rule. At the
center, each block supplies its symbols other than $c$, and the common
cell $(c,c)$ supplies $c$. The global rows and columns remain Latin.

Unequal output symbols still identify their unique block and then their
unique local cell. For equal outputs, the exceptional rule gives
$u+tv=u+t'v$ exactly when $v=0$. Consequently a pair $(a,a)$ occurs at
global cell $(a,c)$ when $a$ belongs to an exceptional block, and at
$(a,a)$ when it belongs to none. For $a=c$, every exceptional block
identifies the same cell $(c,c)$. There is therefore exactly one cell
for every ordered output pair, establishing orthogonality in these
two cases as well.

\subsection{Ray selections}
\label{sec:geo-rays}
A finite direction $d$ defines the ray $\{(r,dr):r\in\F_q^*\}$; the vertical
ray is $\{(0,r):r\in\F_q^*\}$. For disjoint direction sets $P,N,D$, keep
$r$ with $\chi(r)=1$ in directions $P$, $\chi(r)=-1$ in directions $N$,
and every nonzero $r$ in directions $D$. Keep no points on other rays.
Add precisely the infinity points $I_d$ with $d\in I$, and include the origin
only when specified. Every ray recipe uses the ordinary plane.

\paragraph{Order 328 ($q=25$, $k=9$).}
Omit the origin and set
\[
 P=\varnothing,\quad N=\{0,1,2,3,4,\infty\},\quad
 D=\{5,6,\ldots,14\},\quad I=N\cup D.
\]
There are $6\cdot12+10\cdot24+16=328$ points.

\paragraph{Order 335 ($q=25$, $k=9$).}
Let $C=\{d:\operatorname{Nm}(d)=2\}=\{5,12,13,17,18,20\}$.
For $d=a+5b$ put $e_d=1$ when $b\in\{0,1,2\}$ and $e_d=-1$ otherwise;
put $e_\infty=1$. Outside $C$, use a positive or negative half-ray according
to $e_d$. Set $D=\{5,12,13\}$, leave directions $17,18,20$ empty, include
the infinity points of all used rays, and omit the origin. Thus
$n=20\cdot12+3\cdot24+23=335$.

\paragraph{Order 442 ($q=27$, $k=11$).}
Omit the origin and set
\begin{align*}
 P&=\{5,14,16,21\},&N&=\{6,17,22,24,\infty\},\\
 D&=\{0,2,3,4,7,8,9,10,11,12,13,15\},&I&=D\cup\{1\}.
\end{align*}
There are $9\cdot13+12\cdot26+13=442$ points.

\paragraph{Order 494 ($q=31$, $k=11$).}
Include the origin and set
\begin{align*}
 P&=\{2,4,9,12,17,23,25,26,28\},&N&=\{0,1,3\},\\
 D&=\{5,10,11,13,14,15,19,21,22,\infty\},&I&=P\cup N\cup\{\infty\}.
\end{align*}
Then $n=1+12\cdot15+10\cdot30+13=494$.

\paragraph{Ray-intersection counts.}
Let $b_d$ equal 0, 1, or 2 for an empty, half, or full ray, let
$e_d=1,-1,0$ for positive, negative, or other directions, let
$i_d=1$ if $d\in I$ and zero otherwise, and let $o$ be the origin flag.
Put $H=\sum_db_d$ and
\[
 T_d=\sum_{j\in\F_q}\chi(j-d)e_j+e_\infty,
 \qquad T_\infty=\sum_{j\in\F_q}e_j.
\]
The radial line in direction $d$ has $(q-1)b_d/2+o+i_d$ selected points.
The other $q-1$ lines in that direction split equally between the sizes
$(H-b_d+T_d)/2+i_d$ and $(H-b_d-T_d)/2+i_d$.
These formulas follow by solving a line's intersection with each ray;
they give a small integer alternative to the complete line enumeration.

\subsection{Hyperbola selections}
\label{sec:geo-hyperbolas}
For a set $A\subseteq\F_q^*$ select all affine points with $xy\in A$,
then add the listed boundary points. Each nonzero product contributes
$q-1$ points, and all three recipes use the ordinary plane.
\begin{description}
\item[Order 183 ($q=16$, $k=7$).]
Use $A=\{1,2,\ldots,9\}$, both complete axes, and all seventeen infinity
points. The number of points is $9\cdot15+31+17=183$.
\item[Order 430 ($q=27$, $k=11$).]
Use
\[
 A=\{2,3,4,5,6,9,10,11,13,14,15,16,18,19,20,21\}.
\]
Add the origin, no other axis points, and the infinity points of all thirteen
nonsquare finite slopes. The number of points is $16\cdot26+1+13=430$.
\item[Order 469 ($q=27$, $k=11$).]
Use
\[
 A=\{1,5,6,7,8,9,10,11,13,14,15,16,18,19,21,22\}.
\]
Add both complete axes and no infinity points. The count is
$16\cdot26+53=469$.
\end{description}
For an optional count check in the odd-field cases, put
$J_A(t)=\sum_{a\in A}\chi(1+ta)$.
When $m,b\ne0$, the line $y=mx+b$ meets the selected hyperbolas in
$|A|+J_A(4m/b^2)$ points. Add its selected axis or infinity points as specified.
This reduces the general-line calculation to 26 exact character sums;
vertical, horizontal, and zero-intercept lines can be listed separately.

\subsection{Monomial-curve selections}
\label{sec:geo-curves}
A coefficient set $A$ and exponent $e$ specify the points
$(x,cx^e)$ for all $x\ne0$ and $c\in A$. Use field exponentiation, not integer
exponentiation modulo the field size, and then add the following boundary sets.

\paragraph{Order 427 ($q=25$, $k=12$, ordinary plane).}
Put
\[
 E=\{z:\operatorname{Nm}(z)=3\}\cup\{5,20\},\quad
 A=\F_{25}\setminus E,\quad e=4.
\]
Add all $(0,y)$ with $\operatorname{Nm}(y)\ne2$, and no infinity points.
Since $|E|=8$, this gives $24\cdot17+19=427$ points.
Exactly eight blocks have size thirteen. All contain the origin, and any
two intersect only there; fill them by the common-center exception with twelve
squares. All other retained blocks have size at least sixteen and use ordinary
filling.

\paragraph{Order 456 ($q=27$, $k=11$, twisted direction functions).}
Put
\[
 P(c)=c^5-c^3-c^2-c+1,\qquad
 A=\{c\in\F_{27}^*:\chi(P(c))=1\},\qquad e=23.
\]
Equivalently, in integer field labels,
\[
 A=\{2,6,7,8,9,10,11,12,13,14,15,16,17,19,21,24\}.
\]
Add the complete vertical axis and the infinity points of all nonzero square
slopes. There are $16\cdot26+27+13=456$ points. Use ordinary filling.

\paragraph{Order 482 ($q=27$, $k=7$, switched direction functions).}
Let $Q$ be the nonzero squares and
\[
 T=-\{1,\beta^2,(1+\beta)^2\}=\{2,18,23\},\quad A=Q\cup T,\quad e=9.
\]
Add both axes and all thirteen infinity points of nonsquare slopes.
The count is $16\cdot26+53+13=482$.
Exactly thirteen blocks have size eight; they meet pairwise only at the origin.
Use the seven-square common-center filling there, and ordinary filling on the
remaining blocks.

\subsection{Pairs of maximal arcs}
\label{sec:geo-arcs}
A \emph{maximal four-arc} in a plane of order sixteen is a 52-point set meeting
each line in zero or four points. The following quadratic membership tests are
instances of Denniston's construction~\cite{denniston}.
Greig gives the same quadratic-form and additive-subgroup description and
uses complements and nested differences of maximal arcs to construct
PBDs~\cite[Lemma 3.2 and Corollary 3.4]{greig2000}.
Here we use the complement of the symmetric difference of two explicitly
positioned arcs.
Work in the ordinary plane over $\F_{16}$ and set
\[
 A=\{(x,y):x^2+xy+8y^2\in\{0,1,2,3\}\}.
\]
Every integer in this formula is a field label. The set has 52 points and
contains no infinity points.

\paragraph{Order 207 ($k=9$).}
Define a second affine arc
\[
 B=\{(x,y):x^2+x(y+1)+13(y+1)^2\in\{0,6,10,12\}\}.
\]
Keep exactly those ambient points belonging to both arcs or to neither:
$S$ is the complement of the symmetric difference $A\mathbin{\triangle}B$.
The two arcs intersect in nineteen points, so $|S|=273-104+38=207$.
Use ordinary filling with nine squares.

\paragraph{Order 175 ($k=7$).}
For each $(x,y)\in A$, compute
\[
 T(x,y)=\begin{cases}
 (y^{-1},2x/y),&y\ne0,\\ I_{2x},&y=0,
 \end{cases}
\]
where 2 is the field label satisfying $2^4=2+1$. Set $B=T(A)$ and again keep
the complement of $A\mathbin{\triangle}B$ in all 273 ambient points.
This map is the projectivity $[X:Y:Z]\mapsto[Z:2X:Y]$.
The arcs intersect in
$(6,14),(7,3),(12,10)$, giving $|S|=273-104+6=175$.
Use ordinary filling with seven squares.

\subsection{Intersection sizes and verification}
\label{sec:geo-checks}
Table~\ref{tab:geo-histograms} lists the complete intersection histograms.
Notation $s^a$ means that exactly $a$ ambient lines meet the selected set in
$s$ points. The sums of multiplicities are $q^2+q+1$.
The values follow by enumerating the explicit points and lines above. Each
retained size is a prime power supporting the stated local filling; at the two
exceptional orders, additionally check that all minimum blocks contain the
origin. Unique pair incidence then ensures they have no other shared point.
The listed intersection sizes, together with the filling arguments,
establish the twelve bounds.
\Needspace{16\baselineskip}
\begin{center}\small
\begin{longtable}{rrrl}
\caption{Complete line-intersection checks for the geometric selections.}
\label{tab:geo-histograms}\\
\toprule
$n$&$q$&$k$&Histogram\\\midrule\endfirsthead
\toprule $n$&$q$&$k$&Histogram\\\midrule\endhead
175&16&7&$9^{124}\,11^{45}\,13^{101}\,17^3$\\
183&16&7&$9^{45}\,11^{165}\,13^{45}\,17^{18}$\\
207&16&9&$11^{77}\,13^{165}\,17^{31}$\\
328&25&9&$0^{10}\,11^{72}\,13^{486}\,16^{73}\,25^{10}$\\
335&25&9&$0^3\,11^{240}\,13^{164}\,16^{240}\,23^1\,25^3$\\
427&25&12&$0^1\,1^6\,13^8\,16^{246}\,17^{216}\,19^{161}\,25^{13}$\\
430&27&11&$1^2\,13^{183}\,16^{455}\,19^{91}\,25^{26}$\\
442&27&11&$0^6\,1^1\,13^{88}\,16^{325}\,17^{234}\,19^{91}\,27^{12}$\\
456&27&11&$1^1\,13^{157}\,16^{117}\,17^{208}\,19^{247}\,23^{26}\,27^1$\\
469&27&11&$0^1\,13^{78}\,16^{260}\,17^{78}\,19^{312}\,25^{26}\,27^2$\\
482&27&7&$8^{13}\,13^1\,16^{78}\,17^{390}\,19^{234}\,25^{26}\,27^{15}$\\
494&31&11&$1^{10}\,13^{226}\,16^{330}\,17^{327}\,19^{90}\,31^9\,32^1$\\
\bottomrule
\end{longtable}
\end{center}

\clearpage

\appendix
\section{Prompts}\label{app:prompts}

GPT-6 Astra was used in two modes: (a) in ChatGPT Pro with an initial prompt followed by messages to continue, and (b) in the OpenAI Codex harness with the \textit{ultra} reasoning effort, using its \textit{exec} mode for autonomous 12-hour runs. 

The first success came at $n=74$ using the following prompt to Codex:

\begin{Verbatim}[frame=single, breaklines, breaksymbolleft={}, breakindent=0pt]
We are constructing new larger sets of MOLS, improving over existing lower bounds using varied techniques.

Your current task is to construct a set of at least 6 MOLS of order n=74.
Stop when you find this, and output the MOLS in a text file.
Upon successful completion, also provide a brief report summarizing the construction of the result, and identify any code that was used to find the final result.

You can use all 16 cores (32 hyperthreads) of this Linux machine for your work, and you may work on this for up to 12 hours (wall-clock time).

You can use software already installed, including gcc, python3, gap, and standard command-line utilities.  You may install software locally in the working directory.  Install all Python dependencies into a virtual environment under the workspace.  Do not use sudo, apt, pip --user, or system-wide installation.  Build any additional software from source under the workspace and install it under a workspace-local prefix.

Don't get stuck on trying expanding series of repairs to unsuccessful starting points.  Experiment with diverse methods and find a way to successfully solve this instance.  If a starting point isn't working, discard it and find new approaches.
\end{Verbatim}

The subsequent Codex prompts tried to extend this success; the next positive result was 4 MOLS of order $n=22$ using the following prompt:

\begin{Verbatim}[frame=single, breaklines, breaksymbolleft={}, breakindent=0pt]
We are constructing new larger sets of MOLS, improving over existing lower bounds using varied techniques.
We have succeeded in constructing a set of 6 MOLS of order n=74, improving upon the existing public known lower bound for n=74.
The approach used a quasi-difference matrix modulo 67 with quaternion symmetry, completed by filling a seven-symbol hole.
This emerged from a prior effort using Codex, which is in the current directory.  
In particular, the 6 MOLS of order n=74 are in ./output/mols74_6.txt, and the report is in ./output/REPORT.md, and it points to source code used to produce the result.
There are various files in the current directory which are leftover from this previous effort that you won't need, but I preserved the full contents to be sure you have everything you need.

Your current task is to adapt and extend the technique in ./output/REPORT.md to find an improvement in the lower bound for at least one of n in {22,24,26,28,30,33,34,35,36}.
In particular, we're looking for at least the following number of MOLS for each of these n; this is the number of MOLS you need to find, and it would improve upon the existing known lower bound:
n=22: find 4 MOLS
n=24: find 8 MOLS
n=26: find 5 MOLS
n=28: find 6 MOLS
n=30: find 5 MOLS
n=33: find 6 MOLS
n=34: find 5 MOLS
n=35: find 7 MOLS
n=36: find 9 MOLS

Stop when you find at least one of these, and output the MOLS in a text file.
Upon successful completion, also provide a brief report summarizing the construction of the result, and identify any code that was used to find the final result.

You can use all 16 cores (32 hyperthreads) of this Linux machine for your work, and you may work on this for up to 12 hours (wall-clock time).

You can use software already installed, including gcc, python3, gap, and standard command-line utilities.  You can use anything already present in the working directory.  You may install software locally in the working directory.  Install all Python dependencies into a virtual environment under the workspace.  Do not use sudo, apt, pip --user, or system-wide installation.  Build any additional software from source under the workspace and install it under a workspace-local prefix.

Don't get stuck on trying expanding series of repairs to unsuccessful attempts.  Experiment with diverse adaptations and extensions, and find a way to succeed.  If an attempt isn't working, discard it and move on to new adaptations and extensions, or other n.
\end{Verbatim}

After several such improvements with Codex, most of the remaining progress came from long chats with ChatGPT Pro (including the constructions using projective planes).  Below is an example of the initial prompt used.  Note that the ``short program'' idea suggested in the prompt didn't succeed, and the suggested larger improvements in lower bounds didn't materialize, but ChatGPT Pro did produce some of the initial successful constructions based on projective planes in the ensuing long conversation.

\begin{Verbatim}[frame=single, breaklines, breaksymbolleft={}, breakindent=0pt]
Our project is to identify large sets of MOLS (well beyond current known lower bounds) for values of n which are not prime powers.  There is a huge gap between known upper bounds and lower bounds for most values of n that are not prime powers, and our conjecture (which we can assume to be true for purposes of this project) is that lower bounds can be greatly improved for most values of n which are not prime powers.

We've extensively explored constructions based on difference matrices, quasi-difference matrices and separable permutation codes, with some success (using difference matrices and quasi-difference matrices), but we don't think these techniques are really on the right track.  They do succeed in providing a smaller search space, which has enabled slight improvements in some lower bounds for specific n.  But as n grows, the search space still quickly becomes intractable.  It becomes necessary to assume strong further restrictions on the structure of the solution; most such restrictions don't work, and guessing the right structure becomes challenging in itself - and even when there is success, the typical payoff is merely an improvement of 1-4 in the lower bound.  It doesn't seem like that is going to be the ultimate solution to our project.  

Of course, perhaps it could be the case that the search is simply intractable; maybe improved sets of MOLS are going have little structure, and must each be found by an extensive unstructured search which quickly becomes intractable.   But we don't really believe that will be the case.  This isn't like normal NP-hard problems, where (barring P=NP) we expect the time to search for a solution to be exponential in the input size.  Here, the input is just the few bits required to specify n and the number of MOLS we target, so exponential in that isn't large.  It doesn't seem like computational intractability results apply here.

The finite field construction of complete sets of MOLS for prime powers can be viewed as a short program which succeeds quickly for input n that is a prime power.  Our conjecture is that there are other short programs which will succeed quickly (say, in time polynomial in n, and indeed a polynomial that gives practical runtimes for moderate n), for many n which are not prime powers.  It may well be that these programs are are not algebraic in nature at all; they may use other approaches entirely.

Conjecture some approaches we could pursue in this project - what kind of short programs could have this property?  Including non-algebraic approaches?  How could we find these successfully?

We are initially targeting significant improvements in lower bounds for at least a few n<500, for which we have this nice recent table here:
https://www.mdpi.com/2073-8994/16/12/1678
Just to note, we have already improved the lower bounds for n=22, n=60, n=72, n=74, n=225, and n=486 by 1 each; n=192 and n=288 by 2 each; n=162 and n=324 by 3 each, and n=144 by 4 (using difference matrices and quasi-difference matrices).  While this is a nice set of improvements, these were extremely difficult constructions that are unlikely to extend much further.  Also, we are ultimately looking for much larger improvements, that close e.g. 1/3 of the gap or more towards the upper bounds.  This experience was part of what convinced us a radically new approach is needed.
\end{Verbatim}

Conversation continuation messages were mostly simple messages like \texttt{Please continue} and \texttt{Excellent!  Please keep pushing for more lower-bound improvements for n<500, as well as for more ambitious gains.}

\section{The four complete squares of order 22}\label{app:order22}
The four squares below are the fully developed output of
Section~\ref{sec:cyclic22}, with rows and columns in order $0,1,\ldots,21$.
Entries use the alphabet $0,1,\ldots,9,\mathrm{A},\ldots,\mathrm{L}$,
where $\mathrm{A}=10,\ldots,\mathrm{L}=21$; in particular,
$\mathrm{L}$ denotes $\infty$.

\begin{center}
\begin{minipage}[t]{0.485\textwidth}\centering
{\normalsize $L_{0}$}\par\vspace{5pt}
\begingroup\ttfamily\fontsize{10}{11.6}\selectfont
\setlength{\tabcolsep}{1.65pt}\renewcommand{\arraystretch}{1.02}
\begin{tabular}{@{}*{22}{c}@{}}\toprule
% BEGIN_ORDER22_SQUARE_0
L & 4 & D & 9 & 6 & 1 & G & B & 5 & A & 3 & J & H & K & 2 & 7 & E & 8 & C & I & F & 0\\
G & L & 5 & E & A & 7 & 2 & H & C & 6 & B & 4 & K & I & 0 & 3 & 8 & F & 9 & D & J & 1\\
K & H & L & 6 & F & B & 8 & 3 & I & D & 7 & C & 5 & 0 & J & 1 & 4 & 9 & G & A & E & 2\\
F & 0 & I & L & 7 & G & C & 9 & 4 & J & E & 8 & D & 6 & 1 & K & 2 & 5 & A & H & B & 3\\
C & G & 1 & J & L & 8 & H & D & A & 5 & K & F & 9 & E & 7 & 2 & 0 & 3 & 6 & B & I & 4\\
J & D & H & 2 & K & L & 9 & I & E & B & 6 & 0 & G & A & F & 8 & 3 & 1 & 4 & 7 & C & 5\\
D & K & E & I & 3 & 0 & L & A & J & F & C & 7 & 1 & H & B & G & 9 & 4 & 2 & 5 & 8 & 6\\
9 & E & 0 & F & J & 4 & 1 & L & B & K & G & D & 8 & 2 & I & C & H & A & 5 & 3 & 6 & 7\\
7 & A & F & 1 & G & K & 5 & 2 & L & C & 0 & H & E & 9 & 3 & J & D & I & B & 6 & 4 & 8\\
5 & 8 & B & G & 2 & H & 0 & 6 & 3 & L & D & 1 & I & F & A & 4 & K & E & J & C & 7 & 9\\
8 & 6 & 9 & C & H & 3 & I & 1 & 7 & 4 & L & E & 2 & J & G & B & 5 & 0 & F & K & D & A\\
E & 9 & 7 & A & D & I & 4 & J & 2 & 8 & 5 & L & F & 3 & K & H & C & 6 & 1 & G & 0 & B\\
1 & F & A & 8 & B & E & J & 5 & K & 3 & 9 & 6 & L & G & 4 & 0 & I & D & 7 & 2 & H & C\\
I & 2 & G & B & 9 & C & F & K & 6 & 0 & 4 & A & 7 & L & H & 5 & 1 & J & E & 8 & 3 & D\\
4 & J & 3 & H & C & A & D & G & 0 & 7 & 1 & 5 & B & 8 & L & I & 6 & 2 & K & F & 9 & E\\
A & 5 & K & 4 & I & D & B & E & H & 1 & 8 & 2 & 6 & C & 9 & L & J & 7 & 3 & 0 & G & F\\
H & B & 6 & 0 & 5 & J & E & C & F & I & 2 & 9 & 3 & 7 & D & A & L & K & 8 & 4 & 1 & G\\
2 & I & C & 7 & 1 & 6 & K & F & D & G & J & 3 & A & 4 & 8 & E & B & L & 0 & 9 & 5 & H\\
6 & 3 & J & D & 8 & 2 & 7 & 0 & G & E & H & K & 4 & B & 5 & 9 & F & C & L & 1 & A & I\\
B & 7 & 4 & K & E & 9 & 3 & 8 & 1 & H & F & I & 0 & 5 & C & 6 & A & G & D & L & 2 & J\\
3 & C & 8 & 5 & 0 & F & A & 4 & 9 & 2 & I & G & J & 1 & 6 & D & 7 & B & H & E & L & K\\
0 & 1 & 2 & 3 & 4 & 5 & 6 & 7 & 8 & 9 & A & B & C & D & E & F & G & H & I & J & K & L\\
% END_ORDER22_SQUARE_0
\bottomrule\end{tabular}\endgroup
\end{minipage}
\hfill
\begin{minipage}[t]{0.485\textwidth}\centering
{\normalsize $L_{1}$}\par\vspace{5pt}
\begingroup\ttfamily\fontsize{10}{11.6}\selectfont
\setlength{\tabcolsep}{1.65pt}\renewcommand{\arraystretch}{1.02}
\begin{tabular}{@{}*{22}{c}@{}}\toprule
% BEGIN_ORDER22_SQUARE_1
4 & L & H & 8 & G & I & 3 & E & B & F & 0 & 7 & A & C & 9 & 2 & 5 & J & 6 & K & D & 1\\
E & 5 & L & I & 9 & H & J & 4 & F & C & G & 1 & 8 & B & D & A & 3 & 6 & K & 7 & 0 & 2\\
1 & F & 6 & L & J & A & I & K & 5 & G & D & H & 2 & 9 & C & E & B & 4 & 7 & 0 & 8 & 3\\
9 & 2 & G & 7 & L & K & B & J & 0 & 6 & H & E & I & 3 & A & D & F & C & 5 & 8 & 1 & 4\\
2 & A & 3 & H & 8 & L & 0 & C & K & 1 & 7 & I & F & J & 4 & B & E & G & D & 6 & 9 & 5\\
A & 3 & B & 4 & I & 9 & L & 1 & D & 0 & 2 & 8 & J & G & K & 5 & C & F & H & E & 7 & 6\\
8 & B & 4 & C & 5 & J & A & L & 2 & E & 1 & 3 & 9 & K & H & 0 & 6 & D & G & I & F & 7\\
G & 9 & C & 5 & D & 6 & K & B & L & 3 & F & 2 & 4 & A & 0 & I & 1 & 7 & E & H & J & 8\\
K & H & A & D & 6 & E & 7 & 0 & C & L & 4 & G & 3 & 5 & B & 1 & J & 2 & 8 & F & I & 9\\
J & 0 & I & B & E & 7 & F & 8 & 1 & D & L & 5 & H & 4 & 6 & C & 2 & K & 3 & 9 & G & A\\
H & K & 1 & J & C & F & 8 & G & 9 & 2 & E & L & 6 & I & 5 & 7 & D & 3 & 0 & 4 & A & B\\
B & I & 0 & 2 & K & D & G & 9 & H & A & 3 & F & L & 7 & J & 6 & 8 & E & 4 & 1 & 5 & C\\
6 & C & J & 1 & 3 & 0 & E & H & A & I & B & 4 & G & L & 8 & K & 7 & 9 & F & 5 & 2 & D\\
3 & 7 & D & K & 2 & 4 & 1 & F & I & B & J & C & 5 & H & L & 9 & 0 & 8 & A & G & 6 & E\\
7 & 4 & 8 & E & 0 & 3 & 5 & 2 & G & J & C & K & D & 6 & I & L & A & 1 & 9 & B & H & F\\
I & 8 & 5 & 9 & F & 1 & 4 & 6 & 3 & H & K & D & 0 & E & 7 & J & L & B & 2 & A & C & G\\
D & J & 9 & 6 & A & G & 2 & 5 & 7 & 4 & I & 0 & E & 1 & F & 8 & K & L & C & 3 & B & H\\
C & E & K & A & 7 & B & H & 3 & 6 & 8 & 5 & J & 1 & F & 2 & G & 9 & 0 & L & D & 4 & I\\
5 & D & F & 0 & B & 8 & C & I & 4 & 7 & 9 & 6 & K & 2 & G & 3 & H & A & 1 & L & E & J\\
F & 6 & E & G & 1 & C & 9 & D & J & 5 & 8 & A & 7 & 0 & 3 & H & 4 & I & B & 2 & L & K\\
L & G & 7 & F & H & 2 & D & A & E & K & 6 & 9 & B & 8 & 1 & 4 & I & 5 & J & C & 3 & 0\\
0 & 1 & 2 & 3 & 4 & 5 & 6 & 7 & 8 & 9 & A & B & C & D & E & F & G & H & I & J & K & L\\
% END_ORDER22_SQUARE_1
\bottomrule\end{tabular}\endgroup
\end{minipage}
\par\vspace{15pt}
\begin{minipage}[t]{0.485\textwidth}\centering
{\normalsize $L_{2}$}\par\vspace{5pt}
\begingroup\ttfamily\fontsize{10}{11.6}\selectfont
\setlength{\tabcolsep}{1.65pt}\renewcommand{\arraystretch}{1.02}
\begin{tabular}{@{}*{22}{c}@{}}\toprule
% BEGIN_ORDER22_SQUARE_2
G & E & 1 & 5 & C & K & B & 3 & H & D & 7 & I & A & 6 & 4 & F & J & 2 & 9 & 8 & L & 0\\
L & H & F & 2 & 6 & D & 0 & C & 4 & I & E & 8 & J & B & 7 & 5 & G & K & 3 & A & 9 & 1\\
A & L & I & G & 3 & 7 & E & 1 & D & 5 & J & F & 9 & K & C & 8 & 6 & H & 0 & 4 & B & 2\\
C & B & L & J & H & 4 & 8 & F & 2 & E & 6 & K & G & A & 0 & D & 9 & 7 & I & 1 & 5 & 3\\
6 & D & C & L & K & I & 5 & 9 & G & 3 & F & 7 & 0 & H & B & 1 & E & A & 8 & J & 2 & 4\\
3 & 7 & E & D & L & 0 & J & 6 & A & H & 4 & G & 8 & 1 & I & C & 2 & F & B & 9 & K & 5\\
0 & 4 & 8 & F & E & L & 1 & K & 7 & B & I & 5 & H & 9 & 2 & J & D & 3 & G & C & A & 6\\
B & 1 & 5 & 9 & G & F & L & 2 & 0 & 8 & C & J & 6 & I & A & 3 & K & E & 4 & H & D & 7\\
E & C & 2 & 6 & A & H & G & L & 3 & 1 & 9 & D & K & 7 & J & B & 4 & 0 & F & 5 & I & 8\\
J & F & D & 3 & 7 & B & I & H & L & 4 & 2 & A & E & 0 & 8 & K & C & 5 & 1 & G & 6 & 9\\
7 & K & G & E & 4 & 8 & C & J & I & L & 5 & 3 & B & F & 1 & 9 & 0 & D & 6 & 2 & H & A\\
I & 8 & 0 & H & F & 5 & 9 & D & K & J & L & 6 & 4 & C & G & 2 & A & 1 & E & 7 & 3 & B\\
4 & J & 9 & 1 & I & G & 6 & A & E & 0 & K & L & 7 & 5 & D & H & 3 & B & 2 & F & 8 & C\\
9 & 5 & K & A & 2 & J & H & 7 & B & F & 1 & 0 & L & 8 & 6 & E & I & 4 & C & 3 & G & D\\
H & A & 6 & 0 & B & 3 & K & I & 8 & C & G & 2 & 1 & L & 9 & 7 & F & J & 5 & D & 4 & E\\
5 & I & B & 7 & 1 & C & 4 & 0 & J & 9 & D & H & 3 & 2 & L & A & 8 & G & K & 6 & E & F\\
F & 6 & J & C & 8 & 2 & D & 5 & 1 & K & A & E & I & 4 & 3 & L & B & 9 & H & 0 & 7 & G\\
8 & G & 7 & K & D & 9 & 3 & E & 6 & 2 & 0 & B & F & J & 5 & 4 & L & C & A & I & 1 & H\\
2 & 9 & H & 8 & 0 & E & A & 4 & F & 7 & 3 & 1 & C & G & K & 6 & 5 & L & D & B & J & I\\
K & 3 & A & I & 9 & 1 & F & B & 5 & G & 8 & 4 & 2 & D & H & 0 & 7 & 6 & L & E & C & J\\
D & 0 & 4 & B & J & A & 2 & G & C & 6 & H & 9 & 5 & 3 & E & I & 1 & 8 & 7 & L & F & K\\
1 & 2 & 3 & 4 & 5 & 6 & 7 & 8 & 9 & A & B & C & D & E & F & G & H & I & J & K & 0 & L\\
% END_ORDER22_SQUARE_2
\bottomrule\end{tabular}\endgroup
\end{minipage}
\hfill
\begin{minipage}[t]{0.485\textwidth}\centering
{\normalsize $L_{3}$}\par\vspace{5pt}
\begingroup\ttfamily\fontsize{10}{11.6}\selectfont
\setlength{\tabcolsep}{1.65pt}\renewcommand{\arraystretch}{1.02}
\begin{tabular}{@{}*{22}{c}@{}}\toprule
% BEGIN_ORDER22_SQUARE_3
0 & E & I & K & 3 & A & 8 & H & 6 & C & B & 5 & L & J & 2 & 1 & 7 & 4 & F & 9 & D & G\\
E & 1 & F & J & 0 & 4 & B & 9 & I & 7 & D & C & 6 & L & K & 3 & 2 & 8 & 5 & G & A & H\\
B & F & 2 & G & K & 1 & 5 & C & A & J & 8 & E & D & 7 & L & 0 & 4 & 3 & 9 & 6 & H & I\\
I & C & G & 3 & H & 0 & 2 & 6 & D & B & K & 9 & F & E & 8 & L & 1 & 5 & 4 & A & 7 & J\\
8 & J & D & H & 4 & I & 1 & 3 & 7 & E & C & 0 & A & G & F & 9 & L & 2 & 6 & 5 & B & K\\
C & 9 & K & E & I & 5 & J & 2 & 4 & 8 & F & D & 1 & B & H & G & A & L & 3 & 7 & 6 & 0\\
7 & D & A & 0 & F & J & 6 & K & 3 & 5 & 9 & G & E & 2 & C & I & H & B & L & 4 & 8 & 1\\
9 & 8 & E & B & 1 & G & K & 7 & 0 & 4 & 6 & A & H & F & 3 & D & J & I & C & L & 5 & 2\\
6 & A & 9 & F & C & 2 & H & 0 & 8 & 1 & 5 & 7 & B & I & G & 4 & E & K & J & D & L & 3\\
L & 7 & B & A & G & D & 3 & I & 1 & 9 & 2 & 6 & 8 & C & J & H & 5 & F & 0 & K & E & 4\\
F & L & 8 & C & B & H & E & 4 & J & 2 & A & 3 & 7 & 9 & D & K & I & 6 & G & 1 & 0 & 5\\
1 & G & L & 9 & D & C & I & F & 5 & K & 3 & B & 4 & 8 & A & E & 0 & J & 7 & H & 2 & 6\\
3 & 2 & H & L & A & E & D & J & G & 6 & 0 & 4 & C & 5 & 9 & B & F & 1 & K & 8 & I & 7\\
J & 4 & 3 & I & L & B & F & E & K & H & 7 & 1 & 5 & D & 6 & A & C & G & 2 & 0 & 9 & 8\\
A & K & 5 & 4 & J & L & C & G & F & 0 & I & 8 & 2 & 6 & E & 7 & B & D & H & 3 & 1 & 9\\
2 & B & 0 & 6 & 5 & K & L & D & H & G & 1 & J & 9 & 3 & 7 & F & 8 & C & E & I & 4 & A\\
5 & 3 & C & 1 & 7 & 6 & 0 & L & E & I & H & 2 & K & A & 4 & 8 & G & 9 & D & F & J & B\\
K & 6 & 4 & D & 2 & 8 & 7 & 1 & L & F & J & I & 3 & 0 & B & 5 & 9 & H & A & E & G & C\\
H & 0 & 7 & 5 & E & 3 & 9 & 8 & 2 & L & G & K & J & 4 & 1 & C & 6 & A & I & B & F & D\\
G & I & 1 & 8 & 6 & F & 4 & A & 9 & 3 & L & H & 0 & K & 5 & 2 & D & 7 & B & J & C & E\\
D & H & J & 2 & 9 & 7 & G & 5 & B & A & 4 & L & I & 1 & 0 & 6 & 3 & E & 8 & C & K & F\\
4 & 5 & 6 & 7 & 8 & 9 & A & B & C & D & E & F & G & H & I & J & K & 0 & 1 & 2 & 3 & L\\
% END_ORDER22_SQUARE_3
\bottomrule\end{tabular}\endgroup
\end{minipage}
\end{center}
\clearpage
\section{Complete affine-lift certificates}\label{app:affine-certificates}
Every row and column index starts at zero.  An entry in a fine table is
an encoded vector, with base-$p$ digits read least significant first.
The dimensions and decoding rules are those of Section~\ref{sec:affine-lifts}.
\subsection{Order 45}
\begin{center}
\small
\setlength{\tabcolsep}{4pt}
\begin{tabular}{rrrrrrrrr}
\toprule
$t$ & 0 & 1 & 2 & 3 & 4 & 5 & 6 & 7\\
\midrule
0 & 0 & 0 & 0 & 0 & 0 & 0 & 0 & 0\\
1 & 0 & 0 & 0 & 4 & 3 & 2 & 2 & 3\\
2 & 0 & 0 & 3 & 2 & 2 & 4 & 3 & 4\\
3 & 0 & 1 & 2 & 2 & 3 & 3 & 0 & 1\\
4 & 0 & 1 & 3 & 4 & 0 & 1 & 4 & 2\\
5 & 0 & 1 & 4 & 0 & 1 & 4 & 2 & 2\\
6 & 0 & 2 & 1 & 3 & 2 & 3 & 2 & 1\\
7 & 0 & 2 & 2 & 4 & 1 & 0 & 1 & 4\\
8 & 0 & 2 & 4 & 2 & 4 & 1 & 3 & 0\\
9 & 0 & 3 & 1 & 1 & 4 & 4 & 0 & 3\\
10 & 0 & 3 & 2 & 3 & 1 & 2 & 4 & 4\\
11 & 0 & 3 & 4 & 1 & 4 & 3 & 1 & 0\\
12 & 0 & 4 & 0 & 3 & 3 & 1 & 1 & 3\\
13 & 0 & 4 & 1 & 0 & 2 & 0 & 4 & 1\\
14 & 0 & 4 & 3 & 1 & 0 & 2 & 3 & 2\\
\bottomrule
\end{tabular}
\end{center}
\noindent Coarse array $H$ over $\Z_5$.

\begin{center}
\small
\setlength{\tabcolsep}{4pt}
\begin{tabular}{rrrrrrrrr}
\toprule
$t$ & 0 & 1 & 2 & 3 & 4 & 5 & 6 & 7\\
\midrule
0 & 0 & 0 & 0 & 0 & 0 & 0 & 0 & 0\\
1 & 0 & 8 & 7 & 6 & 2 & 3 & 4 & 3\\
2 & 0 & 6 & 2 & 2 & 4 & 6 & 3 & 0\\
3 & 0 & 2 & 2 & 5 & 1 & 0 & 1 & 5\\
4 & 0 & 8 & 5 & 2 & 5 & 7 & 4 & 0\\
5 & 0 & 0 & 1 & 5 & 5 & 2 & 7 & 2\\
6 & 0 & 0 & 4 & 7 & 6 & 3 & 8 & 6\\
7 & 0 & 3 & 4 & 0 & 3 & 8 & 4 & 8\\
8 & 0 & 2 & 7 & 4 & 4 & 0 & 0 & 2\\
9 & 0 & 0 & 3 & 8 & 3 & 4 & 2 & 1\\
10 & 0 & 4 & 1 & 0 & 4 & 7 & 3 & 1\\
11 & 0 & 3 & 4 & 5 & 5 & 8 & 3 & 1\\
12 & 0 & 7 & 6 & 5 & 0 & 3 & 5 & 2\\
13 & 0 & 1 & 5 & 2 & 8 & 3 & 7 & 1\\
14 & 0 & 4 & 0 & 3 & 1 & 0 & 6 & 1\\
\bottomrule
\end{tabular}
\end{center}
\noindent Offsets $c_{t,i}$ encoded in base three.  The one-column map encodings are $(1,3,6,5,4,2,8,7)$.

\Needspace{17\baselineskip}
\subsection{Order 72}
\begin{center}
\small
\setlength{\tabcolsep}{4pt}
\begin{tabular}{rrrrrrrrrr}
\toprule
$t$ & 0 & 1 & 2 & 3 & 4 & 5 & 6 & 7 & 8\\
\midrule
0 & 5 & 3 & 4 & 3 & 3 & 0 & 2 & 5 & 0\\
1 & 0 & 0 & 7 & 7 & 0 & 0 & 7 & 1 & 0\\
2 & 7 & 8 & 5 & 1 & 0 & 0 & 4 & 4 & 7\\
3 & 6 & 2 & 7 & 3 & 6 & 0 & 7 & 2 & 7\\
4 & 8 & 2 & 1 & 1 & 0 & 0 & 4 & 5 & 0\\
5 & 8 & 2 & 4 & 4 & 5 & 0 & 7 & 6 & 5\\
6 & 7 & 0 & 5 & 3 & 2 & 0 & 4 & 5 & 0\\
7 & 7 & 6 & 2 & 7 & 4 & 0 & 4 & 4 & 5\\
\bottomrule
\end{tabular}
\end{center}
\noindent Seed offsets $c_{x,i}$ for layer $t=0$, encoded in base three.  Apply \eqref{eq:72-recurrence} twice.  The coarse permutations and all nine maps are given in the main text.

\subsection{Order 80}
\begin{center}\small
\begin{tabular}{rcc}
\toprule
$s$ & $A_s$ & $B_s$\\
\midrule
0 & $(2,3,2,4,4)$ & $(2,2,3,2,0)$\\
1 & $(0,1,4,4,1)$ & $(3,1,0,1,1)$\\
2 & $(1,0,0,1,3)$ & $(3,0,4,4,2)$\\
3 & $(2,3,3,2,0)$ & $(3,0,1,4,4)$\\
\bottomrule
\end{tabular}
\end{center}
\begin{center}
\small
\setlength{\tabcolsep}{4pt}
\begin{tabular}{rrrrrrrrrrrr}
\toprule
$t$ & 0 & 1 & 2 & 3 & 4 & 5 & 6 & 7 & 8 & 9 & 10\\
\midrule
0 & 0 & 0 & 0 & 0 & 0 & 0 & 0 & 0 & 0 & 0 & 0\\
1 & 0 & 7 & 1 & 2 & 6 & 2 & 2 & 2 & 5 & 3 & 14\\
2 & 0 & 15 & 2 & 7 & 5 & 15 & 1 & 12 & 14 & 10 & 13\\
3 & 0 & 4 & 12 & 12 & 7 & 3 & 8 & 11 & 10 & 1 & 6\\
4 & 0 & 13 & 11 & 15 & 13 & 4 & 5 & 13 & 13 & 6 & 2\\
5 & 0 & 5 & 6 & 9 & 7 & 13 & 11 & 3 & 6 & 0 & 0\\
6 & 0 & 2 & 5 & 15 & 4 & 12 & 9 & 7 & 5 & 3 & 9\\
7 & 0 & 9 & 0 & 4 & 12 & 1 & 15 & 12 & 14 & 7 & 10\\
8 & 0 & 7 & 7 & 8 & 14 & 6 & 8 & 11 & 8 & 12 & 15\\
9 & 0 & 10 & 13 & 11 & 3 & 15 & 5 & 14 & 15 & 2 & 2\\
10 & 0 & 7 & 11 & 13 & 15 & 0 & 6 & 10 & 1 & 6 & 7\\
11 & 0 & 0 & 3 & 3 & 6 & 15 & 13 & 10 & 0 & 14 & 11\\
12 & 0 & 12 & 10 & 7 & 7 & 2 & 13 & 2 & 12 & 4 & 12\\
13 & 0 & 8 & 12 & 15 & 5 & 8 & 1 & 8 & 3 & 9 & 12\\
14 & 0 & 13 & 0 & 12 & 10 & 14 & 14 & 9 & 1 & 9 & 6\\
15 & 0 & 12 & 5 & 1 & 0 & 14 & 12 & 6 & 8 & 15 & 4\\
16 & 0 & 9 & 9 & 2 & 15 & 8 & 7 & 14 & 4 & 5 & 4\\
17 & 0 & 3 & 2 & 3 & 10 & 6 & 11 & 4 & 11 & 5 & 3\\
18 & 0 & 4 & 10 & 13 & 3 & 9 & 4 & 5 & 11 & 8 & 9\\
19 & 0 & 8 & 3 & 9 & 2 & 4 & 12 & 5 & 9 & 7 & 8\\
\bottomrule
\end{tabular}
\end{center}
\noindent Offsets in $\F_{16}$, encoded as four bits.  The slopes are $(0,1,15,10,12,8,11,3,2,13,7)$; multiply by $1$ and $6$ to obtain the two columns of each binary map.

\subsection{Order 96}
\begin{center}
\small
\setlength{\tabcolsep}{4pt}
\begin{tabular}{rrrrrrrrrrrrr}
\toprule
$t$ & 0 & 1 & 2 & 3 & 4 & 5 & 6 & 7 & 8 & 9 & 10 & 11\\
\midrule
0 & 0 & 0 & 0 & 0 & 0 & 0 & 0 & 0 & 0 & 0 & 0 & 0\\
1 & 0 & 0 & 0 & 0 & 1 & 1 & 1 & 1 & 2 & 2 & 2 & 2\\
2 & 0 & 0 & 1 & 2 & 0 & 1 & 2 & 2 & 0 & 1 & 1 & 2\\
3 & 0 & 0 & 2 & 1 & 0 & 2 & 1 & 2 & 1 & 0 & 2 & 1\\
4 & 0 & 1 & 2 & 0 & 2 & 1 & 0 & 2 & 2 & 1 & 0 & 1\\
5 & 0 & 1 & 2 & 1 & 0 & 0 & 2 & 1 & 2 & 2 & 1 & 0\\
6 & 0 & 1 & 0 & 2 & 2 & 2 & 0 & 1 & 1 & 0 & 1 & 2\\
7 & 0 & 1 & 1 & 2 & 2 & 0 & 1 & 0 & 0 & 2 & 2 & 1\\
8 & 0 & 2 & 1 & 0 & 1 & 2 & 2 & 0 & 2 & 0 & 1 & 1\\
9 & 0 & 2 & 1 & 1 & 1 & 0 & 0 & 2 & 1 & 2 & 0 & 2\\
10 & 0 & 2 & 2 & 2 & 1 & 2 & 1 & 1 & 0 & 1 & 0 & 0\\
11 & 0 & 2 & 0 & 1 & 2 & 1 & 2 & 0 & 1 & 1 & 2 & 0\\
\bottomrule
\end{tabular}
\end{center}
\noindent Coarse array $H$ over $\Z_3$.

\begin{center}
\small
\setlength{\tabcolsep}{4pt}
\begin{tabular}{rrrr}
\toprule
$t$ & 0 & 1 & 2\\
\midrule
0 & 2 & 16 & 8\\
1 & 27 & 17 & 2\\
2 & 21 & 4 & 13\\
3 & 15 & 31 & 25\\
4 & 31 & 6 & 23\\
5 & 6 & 8 & 24\\
6 & 19 & 11 & 9\\
7 & 13 & 25 & 4\\
8 & 11 & 2 & 27\\
9 & 4 & 24 & 16\\
10 & 23 & 9 & 6\\
11 & 29 & 21 & 17\\
\bottomrule
\end{tabular}
\end{center}
\noindent Encoded binary columns of the twelve maps $M_i:\F_2^3\to\F_2^5$ (the row index here is $i$).

\begin{center}
\small
\setlength{\tabcolsep}{4pt}
\begin{tabular}{rrrrrrrrrrrrr}
\toprule
$t$ & 0 & 1 & 2 & 3 & 4 & 5 & 6 & 7 & 8 & 9 & 10 & 11\\
\midrule
0 & 9 & 14 & 19 & 28 & 12 & 9 & 23 & 26 & 30 & 9 & 21 & 10\\
1 & 25 & 22 & 26 & 13 & 20 & 25 & 30 & 11 & 15 & 25 & 28 & 18\\
2 & 21 & 22 & 27 & 8 & 23 & 23 & 0 & 22 & 7 & 14 & 19 & 25\\
3 & 13 & 7 & 17 & 11 & 24 & 15 & 14 & 19 & 0 & 30 & 15 & 14\\
4 & 11 & 24 & 10 & 21 & 10 & 26 & 11 & 0 & 13 & 15 & 17 & 3\\
5 & 28 & 12 & 13 & 0 & 5 & 11 & 17 & 15 & 23 & 13 & 12 & 24\\
6 & 12 & 25 & 23 & 5 & 20 & 19 & 31 & 8 & 20 & 21 & 0 & 21\\
7 & 29 & 0 & 30 & 31 & 27 & 31 & 7 & 9 & 30 & 22 & 16 & 11\\
8 & 27 & 29 & 3 & 9 & 13 & 18 & 16 & 26 & 27 & 31 & 26 & 0\\
9 & 19 & 19 & 0 & 18 & 25 & 10 & 21 & 3 & 8 & 23 & 29 & 18\\
10 & 17 & 31 & 10 & 20 & 29 & 17 & 14 & 18 & 22 & 17 & 12 & 27\\
11 & 20 & 15 & 16 & 28 & 0 & 27 & 28 & 29 & 9 & 29 & 5 & 31\\
\bottomrule
\end{tabular}
\end{center}
\noindent Offsets $c_{t,i}\in\F_2^5$, encoded as five bits.

\Needspace{28\baselineskip}
\section{Finite data for the affine array constructions}
\label{app:affine-data}
All indices in these tables start at zero. A displayed digit string is
read from left to right as a vector in the indicated prime field;
it is not an integer in positional notation. In a matrix entry, slashes
separate its rows. Thus \texttt{00100/00010/00001} is a $3$-by-$5$ matrix.
The $C_{t,*}$ entries list the offset vectors in column order. These
conventions apply to all the data used in equation~\eqref{eq:affine-oa}.
{\small
\begin{longtable}{rll}
\caption{Order 54: binary rows and ternary offsets.}\label{tab:affine54-hc}\\
\toprule
$t$ & $H_{t,*}$ & $C_{t,*}$\\
\midrule
\endfirsthead
\toprule
$t$ & $H_{t,*}$ & $C_{t,*}$\\
\midrule
\endhead
0 & \texttt{00000000000} & \texttt{012 012 012 012 012 012 012 012 012 012 012}\\
1 & \texttt{11011100010} & \texttt{110 110 102 100 111 220 202 010 000 000 000}\\
2 & \texttt{11000101101} & \texttt{111 220 202 010 000 000 000 110 110 102 100}\\
3 & \texttt{00101101110} & \texttt{010 000 000 000 110 110 102 100 111 220 202}\\
4 & \texttt{10001011011} & \texttt{220 202 010 000 000 000 110 110 102 100 111}\\
5 & \texttt{10111000101} & \texttt{110 102 100 111 220 202 010 000 000 000 110}\\
6 & \texttt{11100010110} & \texttt{100 111 220 202 010 000 000 000 110 110 102}\\
7 & \texttt{00010110111} & \texttt{202 010 000 000 000 110 110 102 100 111 220}\\
8 & \texttt{01110001011} & \texttt{102 100 111 220 202 010 000 000 000 110 110}\\
9 & \texttt{01101110001} & \texttt{000 110 110 102 100 111 220 202 010 000 000}\\
10 & \texttt{10110111000} & \texttt{000 000 110 110 102 100 111 220 202 010 000}\\
11 & \texttt{01011011100} & \texttt{000 000 000 110 110 102 100 111 220 202 010}\\
\bottomrule
\end{longtable}
}
\Needspace{18\baselineskip}
{\small
\begin{longtable}{rll}
\caption{Order 54: ternary coordinate matrices.}\label{tab:affine54-m}\\
\toprule
$i$ & $A_{i,0}$ & $A_{i,1}$\\
\midrule
\endfirsthead
\toprule
$i$ & $A_{i,0}$ & $A_{i,1}$\\
\midrule
\endhead
0 & \texttt{00100/00010/00001} & \texttt{01100/02010/01001}\\
1 & \texttt{02112/10211/01021} & \texttt{20201/22122/22110}\\
2 & \texttt{20210/22021/12202} & \texttt{22011/20220/11000}\\
3 & \texttt{01201/20120/22012} & \texttt{10212/11112/01020}\\
4 & \texttt{20001/12000/11200} & \texttt{20012/12022/11211}\\
5 & \texttt{21120/22112/02211} & \texttt{02022/11210/10110}\\
6 & \texttt{00122/10012/21001} & \texttt{01012/12122/22221}\\
7 & \texttt{21210/22121/12212} & \texttt{00012/10022/21011}\\
8 & \texttt{22020/22202/02220} & \texttt{12220/02002/22120}\\
9 & \texttt{22102/02210/00221} & \texttt{10122/11220/21211}\\
10 & \texttt{00101/20010/22001} & \texttt{01111/22000/20011}\\
\bottomrule
\end{longtable}
}
\Needspace{34\baselineskip}
{\small
\begin{longtable}{rll}
\caption{Order 98: binary rows and offsets over $\F_7$.}\label{tab:affine98-hc}\\
\toprule
$t$ & $H_{t,*}$ & $C_{t,*}$\\
\midrule
\endfirsthead
\toprule
$t$ & $H_{t,*}$ & $C_{t,*}$\\
\midrule
\endhead
0 & \texttt{000000000} & \texttt{52 31 64 64 52 31 31 64 52}\\
1 & \texttt{111000000} & \texttt{55 36 63 45 16 23 32 61 54}\\
2 & \texttt{000111000} & \texttt{54 32 61 63 55 36 23 45 16}\\
3 & \texttt{000000111} & \texttt{16 23 45 61 54 32 36 63 55}\\
4 & \texttt{111000000} & \texttt{62 51 34 00 00 00 14 22 41}\\
5 & \texttt{000111000} & \texttt{41 14 22 34 62 51 00 00 00}\\
6 & \texttt{000000111} & \texttt{00 00 00 22 41 14 51 34 62}\\
7 & \texttt{100100100} & \texttt{34 46 33 51 25 56 62 13 65}\\
8 & \texttt{010010010} & \texttt{65 62 13 33 34 46 56 51 25}\\
9 & \texttt{001001001} & \texttt{25 56 51 13 65 62 46 33 34}\\
10 & \texttt{011011100} & \texttt{10 33 31 41 33 01 66 03 00}\\
11 & \texttt{101101010} & \texttt{64 20 65 04 14 65 00 53 05}\\
12 & \texttt{001110110} & \texttt{05 00 53 65 64 20 65 04 14}\\
13 & \texttt{110001110} & \texttt{33 01 41 03 00 66 33 31 10}\\
14 & \texttt{110110001} & \texttt{56 52 40 56 02 22 06 00 35}\\
15 & \texttt{101010101} & \texttt{02 22 56 00 35 06 52 40 56}\\
16 & \texttt{010101101} & \texttt{00 66 03 31 10 33 01 41 33}\\
17 & \texttt{011100011} & \texttt{14 65 04 53 05 00 20 65 64}\\
18 & \texttt{100011011} & \texttt{35 06 00 40 56 52 22 56 02}\\
19 & \texttt{101110110} & \texttt{13 02 66 41 22 01 00 00 00}\\
20 & \texttt{110101110} & \texttt{00 00 00 46 04 35 22 14 02}\\
21 & \texttt{011011110} & \texttt{02 22 14 00 00 00 35 46 04}\\
22 & \texttt{110110101} & \texttt{14 41 04 00 00 00 25 01 53}\\
23 & \texttt{011101101} & \texttt{04 35 46 14 02 22 00 00 00}\\
24 & \texttt{101011101} & \texttt{00 00 00 01 53 25 41 04 14}\\
25 & \texttt{011110011} & \texttt{00 00 00 66 13 02 01 41 22}\\
26 & \texttt{101101011} & \texttt{22 01 41 00 00 00 02 66 13}\\
27 & \texttt{110011011} & \texttt{53 25 01 04 14 41 00 00 00}\\
\bottomrule
\end{longtable}
}
{\small
\begin{longtable}{rlll}
\caption{Order 98: slopes over $\F_7$.}\label{tab:affine98-w}\\
\toprule
$i$ & $g=(g_0,g_1)$ & $W_{i,0}$ & $W_{i,1}$\\
\midrule
\endfirsthead
\toprule
$i$ & $g=(g_0,g_1)$ & $W_{i,0}$ & $W_{i,1}$\\
\midrule
\endhead
0 & $(0,0)$ & \texttt{11} & \texttt{43}\\
1 & $(0,1)$ & \texttt{41} & \texttt{23}\\
2 & $(0,2)$ & \texttt{21} & \texttt{13}\\
3 & $(1,0)$ & \texttt{14} & \texttt{45}\\
4 & $(1,1)$ & \texttt{44} & \texttt{25}\\
5 & $(1,2)$ & \texttt{24} & \texttt{15}\\
6 & $(2,0)$ & \texttt{12} & \texttt{46}\\
7 & $(2,1)$ & \texttt{42} & \texttt{26}\\
8 & $(2,2)$ & \texttt{22} & \texttt{16}\\
\bottomrule
\end{longtable}
}

\Needspace{22\baselineskip}
\subsection*{Order 50: compressed binary rows and offsets}
For this table only, the integers in $D$ pack two field elements as
$D=c_0+5c_1$; they are not digit strings in the convention for the
other affine tables.  Apply the blockwise rotation $R$ from the
order-50 construction to the two seed rows simultaneously.
The starting indices and orbit lengths specify the full row order.
{\small
\begin{longtable}{rrll}
\caption{Order 50: complete coarse-row and offset seeds.}
\label{tab:affine50-seeds}\\
\toprule
Start $t$ & Length & $H_{t,*}$ & $D_{t,*}$\\
\midrule
\endfirsthead
\toprule
Start $t$ & Length & $H_{t,*}$ & $D_{t,*}$\\
\midrule
\endhead
0 & 1 & \texttt{000000000} & \texttt{1 1 1 11 11 11 4 4 4}\\
1 & 1 & \texttt{000000000} & \texttt{9 9 9 4 4 4 11 11 11}\\
2 & 3 & \texttt{001001111} & \texttt{18 9 18 17 6 1 9 4 0}\\
5 & 3 & \texttt{001110110} & \texttt{21 4 15 10 3 0 3 10 0}\\
8 & 3 & \texttt{001111001} & \texttt{1 7 7 4 9 0 6 17 1}\\
11 & 3 & \texttt{011011100} & \texttt{0 14 21 6 6 5 24 0 0}\\
14 & 3 & \texttt{011100011} & \texttt{0 4 6 24 0 0 6 5 6}\\
17 & 3 & \texttt{111001001} & \texttt{8 22 0 1 0 0 0 1 0}\\
\bottomrule
\end{longtable}
}

\Needspace{28\baselineskip}
\section{Finite data for the relative difference construction}
\label{app:relative-data}
At order 92 use these fourteen rows, with columns numbered $0$ through
$8$ from left to right. Every row in $H$ is a nine-digit vector over
$\F_7$. Spaces separate the nine entries of each $Y$ row over $\F_{13}$.
Use these rows in the development of Section~\ref{sec:relative-difference}.

\begin{center}
\begin{tabular}{rll}
\toprule
$r$ & $H_{r,*}$ over $\F_7$ & $Y_{r,*}$ over $\F_{13}$\\
\midrule
0  & \texttt{000000000} & \texttt{1 7 4 11 3 9 8 10 12}\\
1  & \texttt{001123356} & \texttt{0 1 5 4 3 6 11 12 8}\\
2  & \texttt{010245654} & \texttt{1 10 8 12 3 9 4 7 11}\\
3  & \texttt{015406526} & \texttt{0 2 7 1 4 5 12 11 3}\\
4  & \texttt{021261433} & \texttt{0 1 10 3 5 9 6 12 8}\\
5  & \texttt{024653240} & \texttt{1 9 10 5 0 2 3 11 4}\\
6  & \texttt{034015163} & \texttt{5 12 10 7 0 9 3 11 2}\\
7  & \texttt{036522405} & \texttt{2 3 7 10 0 12 6 9 8}\\
8  & \texttt{042646335} & \texttt{1 8 7 9 2 0 10 6 4}\\
9  & \texttt{046354121} & \texttt{0 1 12 11 5 10 3 4 9}\\
10 & \texttt{053162014} & \texttt{4 7 12 9 0 2 3 1 10}\\
11 & \texttt{055331642} & \texttt{6 1 3 9 12 0 8 4 10}\\
12 & \texttt{062514512} & \texttt{6 2 7 9 3 0 1 12 4}\\
13 & \texttt{063430261} & \texttt{10 4 1 5 6 0 9 11 8}\\
\bottomrule
\end{tabular}
\end{center}

For order 78 no additional seed matrix is needed:
equations~\eqref{eq:relative78-coarse}--\eqref{eq:relative78-fine}
and the two seven-entry sequences give every seed cell.

\clearpage
\section{Finite data for the additional affine constructions}\label{app:finite}
All vector entries are little-endian prime-base encodings. Indices begin at zero. The row and column orders are those used in Section~\ref{sec:finite-affine}.
\begin{table}[!htbp]
\centering\small
\setlength{\tabcolsep}{3pt}
\renewcommand{\arraystretch}{1.06}
\begin{tabular}{r@{\quad}rrrrrrrrrrrrrrr}
\toprule
$t$ & $0$ & $1$ & $2$ & $3$ & $4$ & $5$ & $6$ & $7$ & $8$ & $9$ & $10$ & $11$ & $12$ & $13$ & $14$\\
\midrule
0 & 0 & 0 & 0 & 0 & 0 & 0 & 0 & 0 & 0 & 0 & 0 & 0 & 0 & 0 & 0\\
1 & 0 & 0 & 0 & 8 & 12 & 14 & 6 & 5 & 7 & 15 & 12 & 13 & 10 & 3 & 15\\
2 & 0 & 8 & 0 & 13 & 11 & 11 & 5 & 9 & 13 & 3 & 14 & 12 & 5 & 4 & 10\\
3 & 0 & 8 & 3 & 11 & 1 & 0 & 10 & 15 & 12 & 2 & 4 & 10 & 5 & 15 & 2\\
4 & 0 & 0 & 13 & 3 & 1 & 15 & 3 & 5 & 12 & 11 & 1 & 1 & 3 & 4 & 14\\
5 & 0 & 0 & 4 & 3 & 4 & 3 & 4 & 9 & 5 & 8 & 2 & 8 & 6 & 15 & 8\\
6 & 0 & 8 & 15 & 9 & 8 & 4 & 5 & 14 & 15 & 6 & 3 & 14 & 7 & 9 & 9\\
7 & 0 & 8 & 11 & 14 & 9 & 3 & 1 & 8 & 8 & 12 & 1 & 14 & 1 & 3 & 14\\
8 & 0 & 8 & 12 & 3 & 3 & 1 & 15 & 7 & 12 & 12 & 12 & 6 & 2 & 6 & 9\\
9 & 0 & 8 & 3 & 12 & 8 & 13 & 13 & 4 & 11 & 7 & 12 & 6 & 11 & 10 & 13\\
10 & 0 & 8 & 1 & 2 & 8 & 1 & 1 & 12 & 10 & 12 & 13 & 11 & 15 & 5 & 0\\
11 & 0 & 0 & 11 & 14 & 2 & 1 & 3 & 3 & 13 & 10 & 12 & 11 & 4 & 10 & 5\\
12 & 0 & 0 & 15 & 11 & 13 & 2 & 14 & 7 & 2 & 11 & 12 & 8 & 11 & 5 & 9\\
13 & 0 & 0 & 0 & 5 & 5 & 13 & 7 & 12 & 0 & 13 & 3 & 14 & 13 & 7 & 15\\
14 & 0 & 0 & 4 & 12 & 6 & 4 & 7 & 2 & 3 & 5 & 7 & 9 & 5 & 0 & 7\\
15 & 0 & 8 & 7 & 1 & 15 & 2 & 4 & 6 & 15 & 4 & 13 & 13 & 3 & 2 & 5\\
16 & 0 & 0 & 6 & 14 & 13 & 0 & 13 & 11 & 3 & 0 & 11 & 15 & 12 & 13 & 6\\
17 & 0 & 8 & 2 & 12 & 9 & 5 & 8 & 12 & 9 & 9 & 7 & 8 & 13 & 5 & 12\\
\bottomrule
\end{tabular}
\caption{The 18 offset rows for order 144; entries encode four binary digits.}
\label{tab:finite-c144}
\end{table}
\begin{table}[!htbp]
\centering\small
\setlength{\tabcolsep}{3pt}
\renewcommand{\arraystretch}{1.06}
\begin{tabular}{r@{\quad}rrr}
\toprule
$i$ & $0$ & $1$ & $2$\\
\midrule
0 & 45 & 61 & 26\\
1 & 60 & 23 & 69\\
2 & 5 & 15 & 45\\
3 & 23 & 69 & 50\\
4 & 37 & 28 & 1\\
5 & 19 & 57 & 14\\
6 & 13 & 39 & 43\\
7 & 1 & 3 & 9\\
8 & 65 & 29 & 4\\
9 & 3 & 9 & 27\\
10 & 52 & 73 & 62\\
\bottomrule
\end{tabular}
\caption{The encoded columns of the eleven $4\times3$ matrices for order 162. Decode each integer into four ternary digits.}
\label{tab:finite-m162}
\end{table}
\begin{table}[!htbp]
\centering\small
\setlength{\tabcolsep}{3pt}
\renewcommand{\arraystretch}{1.06}
\begin{tabular}{r@{\quad}rrrrrrrrrrr}
\toprule
$t$ & $0$ & $1$ & $2$ & $3$ & $4$ & $5$ & $6$ & $7$ & $8$ & $9$ & $10$\\
\midrule
0 & 47 & 24 & 39 & 34 & 16 & 23 & 46 & 42 & 48 & 0 & 0\\
1 & 80 & 61 & 14 & 55 & 62 & 37 & 20 & 9 & 43 & 1 & 0\\
2 & 48 & 50 & 55 & 51 & 67 & 64 & 35 & 7 & 63 & 2 & 0\\
3 & 44 & 46 & 77 & 52 & 22 & 79 & 64 & 46 & 72 & 2 & 0\\
4 & 20 & 71 & 23 & 50 & 65 & 68 & 33 & 38 & 73 & 2 & 0\\
5 & 67 & 69 & 27 & 19 & 23 & 54 & 60 & 23 & 10 & 1 & 0\\
6 & 49 & 77 & 43 & 69 & 76 & 52 & 21 & 28 & 16 & 0 & 0\\
7 & 8 & 30 & 66 & 50 & 11 & 60 & 2 & 42 & 33 & 1 & 0\\
8 & 0 & 0 & 50 & 47 & 72 & 27 & 76 & 44 & 59 & 2 & 0\\
9 & 64 & 73 & 0 & 80 & 57 & 58 & 0 & 0 & 48 & 1 & 0\\
10 & 74 & 57 & 77 & 0 & 0 & 0 & 3 & 0 & 0 & 0 & 0\\
11 & 0 & 0 & 0 & 0 & 0 & 0 & 0 & 0 & 0 & 0 & 0\\
\bottomrule
\end{tabular}
\caption{The 12 offset rows for order 162; entries encode four ternary digits.}
\label{tab:finite-c162}
\end{table}
\begin{table}[!htbp]
\centering\small
\setlength{\tabcolsep}{3pt}
\renewcommand{\arraystretch}{1.06}
\begin{tabular}{r@{\quad}rrrrrrrrrr}
\toprule
$t$ & $0$ & $1$ & $2$ & $3$ & $4$ & $5$ & $6$ & $7$ & $8$ & $9$\\
\midrule
0 & 0 & 0 & 0 & 0 & 0 & 0 & 0 & 0 & 0 & 0\\
1 & 0 & 0 & 0 & 0 & 1 & 1 & 1 & 1 & 2 & 2\\
2 & 0 & 0 & 1 & 2 & 0 & 1 & 2 & 2 & 0 & 1\\
3 & 0 & 0 & 2 & 1 & 0 & 2 & 1 & 2 & 1 & 0\\
4 & 0 & 1 & 2 & 0 & 2 & 1 & 0 & 2 & 2 & 1\\
5 & 0 & 1 & 2 & 1 & 0 & 0 & 2 & 1 & 2 & 2\\
6 & 0 & 1 & 0 & 2 & 2 & 2 & 0 & 1 & 1 & 0\\
7 & 0 & 1 & 1 & 2 & 2 & 0 & 1 & 0 & 0 & 2\\
8 & 0 & 2 & 1 & 0 & 1 & 2 & 2 & 0 & 2 & 0\\
9 & 0 & 2 & 1 & 1 & 1 & 0 & 0 & 2 & 1 & 2\\
10 & 0 & 2 & 2 & 2 & 1 & 2 & 1 & 1 & 0 & 1\\
11 & 0 & 2 & 0 & 1 & 2 & 1 & 2 & 0 & 1 & 1\\
\bottomrule
\end{tabular}
\caption{The 12 coarse rows for order 192, over $\Z_3$.}
\label{tab:finite-h192}
\end{table}
\begin{table}[!htbp]
\centering\small
\setlength{\tabcolsep}{3pt}
\renewcommand{\arraystretch}{1.06}
\begin{tabular}{r@{\quad}rrrrrrrrrr}
\toprule
$t$ & $0$ & $1$ & $2$ & $3$ & $4$ & $5$ & $6$ & $7$ & $8$ & $9$\\
\midrule
0 & 0 & 31 & 26 & 48 & 8 & 39 & 45 & 54 & 32 & 0\\
1 & 10 & 32 & 23 & 14 & 38 & 20 & 61 & 28 & 16 & 0\\
2 & 36 & 37 & 52 & 22 & 59 & 10 & 29 & 38 & 16 & 0\\
3 & 10 & 53 & 0 & 19 & 4 & 14 & 42 & 52 & 16 & 0\\
4 & 62 & 35 & 9 & 27 & 21 & 33 & 62 & 6 & 0 & 0\\
5 & 63 & 2 & 15 & 3 & 30 & 14 & 11 & 6 & 0 & 0\\
6 & 62 & 31 & 10 & 38 & 0 & 2 & 0 & 32 & 0 & 32\\
7 & 1 & 4 & 0 & 0 & 48 & 49 & 1 & 16 & 32 & 20\\
8 & 1 & 0 & 32 & 25 & 3 & 19 & 20 & 22 & 49 & 28\\
9 & 0 & 32 & 56 & 55 & 25 & 50 & 19 & 50 & 23 & 41\\
10 & 0 & 48 & 29 & 36 & 59 & 25 & 53 & 56 & 45 & 34\\
11 & 0 & 16 & 49 & 0 & 60 & 42 & 12 & 55 & 56 & 14\\
\bottomrule
\end{tabular}
\caption{The 12 offset rows for order 192; entries encode six binary digits.}
\label{tab:finite-c192}
\end{table}
\begin{table}[!htbp]
\centering\small
\setlength{\tabcolsep}{3pt}
\renewcommand{\arraystretch}{1.06}
\begin{tabular}{r@{\quad}rrrrrrrrrrrrrrrrrr}
\toprule
$t$ & $0$ & $1$ & $2$ & $3$ & $4$ & $5$ & $6$ & $7$ & $8$ & $9$ & $10$ & $11$ & $12$ & $13$ & $14$ & $15$ & $16$ & $17$\\
\midrule
0 & 0 & 0 & 0 & 0 & 0 & 0 & 0 & 0 & 0 & 0 & 0 & 0 & 0 & 0 & 0 & 0 & 0 & 0\\
1 & 0 & 0 & 0 & 0 & 0 & 0 & 0 & 0 & 0 & 16 & 1 & 4 & 16 & 8 & 4 & 4 & 16 & 16\\
2 & 0 & 16 & 4 & 5 & 8 & 27 & 21 & 4 & 8 & 28 & 24 & 4 & 14 & 11 & 25 & 10 & 27 & 12\\
3 & 0 & 16 & 5 & 24 & 29 & 15 & 15 & 19 & 6 & 21 & 7 & 31 & 17 & 19 & 16 & 3 & 3 & 21\\
4 & 0 & 0 & 0 & 30 & 24 & 4 & 27 & 18 & 27 & 30 & 22 & 15 & 14 & 29 & 5 & 2 & 24 & 29\\
5 & 0 & 16 & 9 & 28 & 7 & 17 & 10 & 22 & 7 & 30 & 19 & 12 & 31 & 31 & 19 & 9 & 6 & 22\\
6 & 0 & 16 & 26 & 26 & 16 & 14 & 11 & 14 & 1 & 17 & 27 & 24 & 6 & 14 & 6 & 26 & 14 & 27\\
7 & 0 & 0 & 1 & 18 & 23 & 22 & 15 & 30 & 1 & 9 & 19 & 12 & 2 & 18 & 12 & 15 & 31 & 15\\
8 & 0 & 16 & 26 & 20 & 20 & 6 & 1 & 30 & 7 & 2 & 16 & 26 & 14 & 17 & 17 & 31 & 18 & 22\\
9 & 0 & 16 & 15 & 1 & 18 & 29 & 15 & 26 & 20 & 31 & 4 & 29 & 9 & 6 & 13 & 19 & 25 & 13\\
10 & 0 & 16 & 29 & 17 & 3 & 15 & 8 & 11 & 16 & 26 & 9 & 9 & 31 & 20 & 8 & 17 & 20 & 26\\
11 & 0 & 0 & 31 & 25 & 9 & 23 & 21 & 26 & 23 & 7 & 19 & 9 & 12 & 10 & 25 & 6 & 23 & 9\\
12 & 0 & 16 & 23 & 11 & 18 & 26 & 19 & 13 & 9 & 28 & 0 & 1 & 8 & 13 & 27 & 26 & 25 & 29\\
13 & 0 & 0 & 28 & 11 & 11 & 12 & 9 & 17 & 16 & 9 & 22 & 10 & 26 & 9 & 12 & 20 & 30 & 28\\
14 & 0 & 0 & 15 & 22 & 6 & 14 & 30 & 19 & 23 & 20 & 2 & 3 & 27 & 7 & 0 & 2 & 28 & 11\\
15 & 0 & 16 & 26 & 15 & 24 & 0 & 28 & 25 & 28 & 2 & 8 & 11 & 28 & 7 & 6 & 25 & 14 & 0\\
16 & 0 & 0 & 3 & 1 & 22 & 29 & 5 & 17 & 28 & 26 & 5 & 17 & 7 & 20 & 19 & 17 & 20 & 2\\
17 & 0 & 0 & 28 & 7 & 26 & 5 & 14 & 14 & 28 & 28 & 21 & 28 & 2 & 25 & 5 & 7 & 3 & 11\\
\bottomrule
\end{tabular}
\caption{The 18 offset rows for order 288; entries encode five binary digits.}
\label{tab:finite-c288}
\end{table}
\begin{table}[!htbp]
\centering\small
\setlength{\tabcolsep}{3pt}
\renewcommand{\arraystretch}{1.06}
\begin{tabular}{r@{\quad}rrrrrrrrrrrr}
\toprule
$t$ & $0$ & $1$ & $2$ & $3$ & $4$ & $5$ & $6$ & $7$ & $8$ & $9$ & $10$ & $11$\\
\midrule
0 & 0 & 0 & 0 & 1 & 1 & 1 & 2 & 2 & 2 & 3 & 3 & 3\\
1 & 0 & 0 & 0 & 2 & 2 & 2 & 3 & 3 & 3 & 1 & 1 & 1\\
2 & 0 & 0 & 0 & 3 & 3 & 3 & 1 & 1 & 1 & 2 & 2 & 2\\
3 & 1 & 2 & 3 & 1 & 2 & 3 & 1 & 2 & 3 & 1 & 2 & 3\\
4 & 1 & 2 & 3 & 2 & 3 & 1 & 2 & 3 & 1 & 2 & 3 & 1\\
5 & 1 & 2 & 3 & 3 & 1 & 2 & 3 & 1 & 2 & 3 & 1 & 2\\
6 & 2 & 3 & 1 & 1 & 2 & 3 & 2 & 3 & 1 & 3 & 1 & 2\\
7 & 2 & 3 & 1 & 2 & 3 & 1 & 3 & 1 & 2 & 1 & 2 & 3\\
8 & 2 & 3 & 1 & 3 & 1 & 2 & 1 & 2 & 3 & 2 & 3 & 1\\
9 & 3 & 1 & 2 & 1 & 2 & 3 & 3 & 1 & 2 & 2 & 3 & 1\\
10 & 3 & 1 & 2 & 2 & 3 & 1 & 1 & 2 & 3 & 3 & 1 & 2\\
11 & 3 & 1 & 2 & 3 & 1 & 2 & 2 & 3 & 1 & 1 & 2 & 3\\
\bottomrule
\end{tabular}
\caption{The 12 coarse rows for order 324, over $\Z_2^2$ (two-bit XOR).}
\label{tab:finite-h324}
\end{table}
\begin{table}[!htbp]
\centering\small
\setlength{\tabcolsep}{3pt}
\renewcommand{\arraystretch}{1.06}
\begin{tabular}{r@{\quad}rrrrrrrrrrrr}
\toprule
$t$ & $0$ & $1$ & $2$ & $3$ & $4$ & $5$ & $6$ & $7$ & $8$ & $9$ & $10$ & $11$\\
\midrule
0 & 59 & 56 & 68 & 73 & 6 & 1 & 56 & 48 & 5 & 30 & 27 & 0\\
1 & 55 & 27 & 22 & 4 & 58 & 21 & 65 & 61 & 60 & 66 & 0 & 0\\
2 & 56 & 21 & 51 & 80 & 58 & 27 & 28 & 57 & 24 & 0 & 54 & 0\\
3 & 61 & 75 & 18 & 29 & 74 & 1 & 64 & 54 & 34 & 48 & 54 & 0\\
4 & 50 & 57 & 53 & 19 & 32 & 2 & 74 & 29 & 4 & 24 & 54 & 0\\
5 & 74 & 53 & 44 & 3 & 47 & 71 & 66 & 25 & 6 & 12 & 54 & 0\\
6 & 43 & 36 & 40 & 77 & 54 & 68 & 29 & 74 & 19 & 48 & 27 & 0\\
7 & 76 & 24 & 32 & 33 & 4 & 4 & 32 & 20 & 48 & 78 & 27 & 0\\
8 & 49 & 15 & 66 & 65 & 45 & 2 & 0 & 52 & 43 & 60 & 0 & 0\\
9 & 38 & 76 & 1 & 16 & 69 & 41 & 67 & 65 & 42 & 48 & 27 & 0\\
10 & 75 & 34 & 4 & 72 & 75 & 11 & 2 & 67 & 9 & 54 & 0 & 0\\
11 & 0 & 0 & 0 & 0 & 0 & 0 & 0 & 0 & 0 & 0 & 0 & 0\\
\bottomrule
\end{tabular}
\caption{The 12 offset rows for order 324; entries encode four ternary digits.}
\label{tab:finite-c324}
\end{table}
\begin{table}[!htbp]
\centering\small
\setlength{\tabcolsep}{3pt}
\renewcommand{\arraystretch}{1.06}
\begin{tabular}{r@{\quad}rrrrrrrrrrr}
\toprule
$t$ & $0$ & $1$ & $2$ & $3$ & $4$ & $5$ & $6$ & $7$ & $8$ & $9$ & $10$\\
\midrule
0 & 0 & 0 & 0 & 0 & 0 & 0 & 0 & 0 & 0 & 0 & 0\\
1 & 1 & 0 & 1 & 0 & 0 & 0 & 1 & 1 & 1 & 0 & 1\\
2 & 1 & 1 & 0 & 1 & 0 & 0 & 0 & 1 & 1 & 1 & 0\\
3 & 0 & 1 & 1 & 0 & 1 & 0 & 0 & 0 & 1 & 1 & 1\\
4 & 1 & 0 & 1 & 1 & 0 & 1 & 0 & 0 & 0 & 1 & 1\\
5 & 1 & 1 & 0 & 1 & 1 & 0 & 1 & 0 & 0 & 0 & 1\\
6 & 1 & 1 & 1 & 0 & 1 & 1 & 0 & 1 & 0 & 0 & 0\\
7 & 0 & 1 & 1 & 1 & 0 & 1 & 1 & 0 & 1 & 0 & 0\\
8 & 0 & 0 & 1 & 1 & 1 & 0 & 1 & 1 & 0 & 1 & 0\\
9 & 0 & 0 & 0 & 1 & 1 & 1 & 0 & 1 & 1 & 0 & 1\\
10 & 1 & 0 & 0 & 0 & 1 & 1 & 1 & 0 & 1 & 1 & 0\\
11 & 0 & 1 & 0 & 0 & 0 & 1 & 1 & 1 & 0 & 1 & 1\\
\bottomrule
\end{tabular}
\caption{The common binary array $B$ for orders 162 and 486. Use the column permutation stated in the text at order 162.}
\label{tab:finite-binary}
\end{table}
\begin{table}[!htbp]
\centering\small
\setlength{\tabcolsep}{3pt}
\renewcommand{\arraystretch}{1.06}
\begin{tabular}{r@{\quad}rrrrrrrrrrr}
\toprule
$t$ & $0$ & $1$ & $2$ & $3$ & $4$ & $5$ & $6$ & $7$ & $8$ & $9$ & $10$\\
\midrule
0 & 216 & 43 & 100 & 155 & 9 & 178 & 3 & 124 & 218 & 2 & 0\\
1 & 173 & 189 & 166 & 149 & 113 & 225 & 162 & 28 & 188 & 2 & 0\\
2 & 211 & 40 & 61 & 53 & 180 & 167 & 220 & 102 & 140 & 2 & 0\\
3 & 148 & 34 & 77 & 140 & 81 & 36 & 181 & 69 & 165 & 2 & 0\\
4 & 102 & 171 & 182 & 5 & 13 & 15 & 20 & 239 & 187 & 1 & 0\\
5 & 81 & 17 & 123 & 204 & 200 & 90 & 229 & 72 & 157 & 1 & 0\\
6 & 106 & 51 & 226 & 242 & 138 & 83 & 47 & 185 & 75 & 1 & 0\\
7 & 195 & 187 & 128 & 86 & 188 & 107 & 56 & 213 & 195 & 0 & 0\\
8 & 22 & 200 & 0 & 140 & 78 & 0 & 154 & 82 & 16 & 1 & 0\\
9 & 21 & 85 & 197 & 0 & 8 & 23 & 0 & 2 & 0 & 0 & 0\\
10 & 0 & 0 & 102 & 14 & 0 & 3 & 2 & 0 & 0 & 0 & 0\\
11 & 0 & 0 & 0 & 0 & 0 & 0 & 0 & 0 & 0 & 0 & 0\\
\bottomrule
\end{tabular}
\caption{The 12 offset rows for order 486; entries encode five ternary digits.}
\label{tab:finite-c486}
\end{table}

\end{document}